\documentclass[11pt,francais]{article}
\usepackage{amsmath}\usepackage{epsf,amsfonts,amsthm}\usepackage{amscd,mathptmx,amssymb}
\usepackage{xcolor,epic,eepic}\usepackage{epsfig}
\usepackage{fontenc,indentfirst, delarray,amsfonts,amsmath,amssymb}
\usepackage{rotating}
\usepackage[T1]{fontenc}

\newcommand{\bee}{\begin{enumerate}}
\newcommand{\eee}{\end{enumerate}}
\newcommand{\benn}{\begin{equation*}}
\newcommand{\eenn}{\end{equation*}}
\newcommand{\be}{\begin{equation}}
\newcommand{\ee}{\end{equation}}
\newcommand{\bean}{\begin{eqnarray}}
\newcommand{\eean}{\end{eqnarray}}
\newcommand{\bea}{\begin{eqnarray*}}
\newcommand{\eea}{\end{eqnarray*}}

\newcommand{\ot}{\otimes}

\newcommand{\p}{\partial}
\newcommand{\la}{\langle}
\newcommand{\ra}{\rangle}
\newcommand{\raa}{\rightarrow}
\newcommand{\lraa}{\longrightarrow}
\newcommand{\Ci}{C^{\infty}}
\newcommand{\E}{\ell}
\newcommand{\N}{\mathbb{N}}
\newcommand{\Z}{\mathbb{Z}}

\newcommand{\K}{\mathbb{K}}

\newcommand{\cf}{{\cal F}}

\newcommand{\lp}{\left(}
\newcommand{\rp}{\right)}

\newcommand{\A}{{\cal A}}

\newcommand{\op}[1]{\!\!\mathop{\rm ~#1}\nolimits}

\newcommand{\da}{\downarrow}
\newcommand{\ua}{\uparrow}

\newcommand{\emp}{\emptyset}

\newcommand{\cg}{{\cal G}}

\mathchardef\za="710B  %\alpha
\mathchardef\zb="710C  %\beta
\mathchardef\zg="710D  %\gamma
\mathchardef\zd="710E  %\delta
\mathchardef\zve="710F %\epsilon
\mathchardef\zz="7110  %\zeta
\mathchardef\zh="7111  %\eta
\mathchardef\zy="7112 %\theta

\mathchardef\zi="7113  %\iota
\mathchardef\zk="7114  %\kappa
\mathchardef\zl="7115  %\lambda
\mathchardef\zm="7116  %\mu
\mathchardef\zn="7117  %\nu
\mathchardef\zx="7118  %\xi
\mathchardef\zp="7119  %\pi
\mathchardef\zr="711A  %\rho
\mathchardef\zs="711B  %\sigma
\mathchardef\zt="711C  %\tau
\mathchardef\zu="711D  %\upsilon
\mathchardef\zf="711E %\phi
\mathchardef\zq="711F  %\chi
\mathchardef\zc="7120  %\psi
\mathchardef\zw="7121  %\omega
\mathchardef\ze="7122  %\varepsilon
\mathchardef\zvy="7123  %\vartheta
\mathchardef\zvw="7124  %\varomega
\mathchardef\zvr="7125 %\varrho
\mathchardef\zvs="7126 %\varsigma
\mathchardef\zvf="7127  %\varphi
\mathchardef\zG="7000  %\Gamma
\mathchardef\zD="7001  %\Delta
\mathchardef\zY="7002  %\Theta
\mathchardef\zL="7003  %\Lambda
\mathchardef\zX="7004  %\Xi
\mathchardef\zP="7005  %\Pi
\mathchardef\zS="7006  %\Sigma
\mathchardef\zU="7007  %\Upsilon
\mathchardef\zF="7008  %\Phi
\mathchardef\zW="700A  %\Omega

\newcommand{\cyclic}{\mathop{\kern0.9ex{{+}
 \kern-2.15ex\raise-.25ex\hbox{\Large\hbox{$\circlearrowright$}}}}\limits}
\newtheorem{rem}{Remark}
\newtheorem{theo}{Theorem}
\newtheorem{prop}{Proposition}
\newtheorem{lem}{Lemma}
\newtheorem{cor}{Corollary}
\newtheorem{ex}{Example}
\newtheorem{defi}{Definition}

\begin{document}
\title{Coalgebraic Approach to the Loday Infinity Category,\\ Stem Differential for $2n$-ary Graded and Homotopy Algebras}
\author{Mourad Ammar, Norbert Poncin \footnote{University of Luxembourg, Campus Limpertsberg,
Mathematics Research Unit, 162A, avenue de la Fa\"iencerie, L-1511
Luxembourg City, Grand-Duchy of Luxembourg, E-mail:
mourad.ammar@uni.lu, norbert.poncin@uni.lu. The research of N.
Poncin was supported by UL-grant SGQnp2008.}}\maketitle

\begin{abstract}
\noindent We define a graded twisted-coassociative coproduct on
the tensor algebra $TW$ of any $\Z^n$-graded vector space $W$. If
$W$ is the desuspension space $\da V$ of a graded vector space
$V$, the coderivations (resp. quadratic ``degree $1$''
codifferentials, arbitrary odd codifferentials) of this coalgebra
are $1$-to-$1$ with sequences $\zp_s$, $s\ge 1$, of $s$-linear
maps on $V$ (resp. $\Z^n$-graded Loday structures on $V$,
sequences that we call Loday infinity structures on $V$). We prove
a minimal model theorem for Loday infinity algebras, investigate
Loday infinity morphisms, and observe that the $\op{Lod}_{\infty}$
category contains the $\op{L}_{\infty}$ category as a subcategory.
Moreover, the graded Lie bracket of coderivations gives rise to a
graded Lie ``stem'' bracket on the cochain spaces of graded Loday,
Loday infinity, and $2n$-ary graded Loday algebras (the latter
extend the corresponding Lie algebras in the sense of Michor and
Vinogradov). These algebraic structures have square zero with
respect to the stem bracket, so that we obtain natural
cohomological theories that have good properties with respect to
formal deformations. The stem bracket restricts to the graded
Nijenhuis-Richardson and---up to isomorphism---to the
Grabowski-Marmo brackets (the last bracket extends the
Schouten-Nijenhuis bracket to the space of graded antisymmetric
first order polydifferential operators), and it encodes, beyond
the already mentioned cohomologies, those of graded Lie, graded
Poisson, graded Jacobi, Lie infinity, as well as that of $2n$-ary
graded Lie algebras.\end{abstract}

\vspace{2mm}


\section{Introduction}

The prevailing approach to cohomology consists in associating in a
canonical way a differential module to the investigated algebraic
structure. For instance, in \cite{LLMP2} and \cite{ILLMP} (resp.
in \cite{Oud}) the author(s) define(s) Lichnerowicz-Jacobi and
Nambu-Poisson cohomologies (resp. the cohomology of graded Leibniz
algebras with trivial coefficients) essentially as the cohomology
of a Lie algebroid (resp. as the homology of a differential graded
dual Leibniz algebra) that is associated with any Jacobi or
Nambu-Poisson manifold [it suffices to solve the Lie algebra
homomorphism equation for the natural anchor map] (resp. that is
induced by the considered Leibniz algebra [via
dualization]).\medskip

However, the most natural way to obtain a differential space from
an algebraic structure is the square-zero or canonical element
method that was initiated by De Wilde and Lecomte in \cite{DWL}.
They look for a graded Lie algebra, such that there is a
$1$-to-$1$ correspondence between the studied algebraic structures
and the degree $1$ elements of the Lie algebra space that square
to zero with respect to the Lie bracket. The adjoint action by the
given structure then provides a differential graded Lie algebra,
the cohomology of which has the usual properties with respect to
formal deformations of the algebraic structure.\medskip

Two efficient tools allow finding the mentioned graded Lie
algebra:\medskip

- The coalgebraic technique consists in the identification of the
cochain space of the investigated algebra with certain
coderivations, in such a way that the algebraic structure can be
viewed as an odd codifferential. These identifications enable
constructing a graded Lie bracket on the space of cochains by
transfer of the commutator bracket of coderivations; the
considered algebraic structures are then canonical elements of
this transferred bracket. The procedure was applied by Stasheff
\cite{St} in the case of associative and Lie algebras and by
Penkava \cite{Pen} for strongly homotopy associative and Lie
algebras.

- The operadic theory, as developed in \cite{MSS}, shows that
Stasheff's approach is applicable for any quadratic operad $\cal
P$. The authors associate to ${\cal P}$ a cofree coalgebra over
the dual operad, whose quadratic codifferentials correspond to the
${\cal P}-$algebra structures on a graded vector space $V$. This
construction yields the homology and cohomology theories of ${\cal
P}-$algebras on $V$ and allows defining strongly homotopy ${\cal
P}-$algebra structures on $V$ as arbitrary
codifferentials.\medskip
%The operadic method was used by Balavoine \cite{Bal2}, who
%constructed, for any quadratic operad ${\cal P}$ and any
%finite-dimensional vector space $V$, by means of the Koszul dual
%in the sense of Ginzburg and Kapranov \cite{GK}, a graded Lie
%algebra, whose canonical elements are $1$-to-$1$ with the ${\cal
%P}$-algebra structures on $V$. He then applied his theory to
%nongraded associative, commutative, Lie, and Loday
%algebras.\medskip

In the present work we provide a tensor coalgebra that induces the
proper concepts of Loday infinity algebras and morphisms, and
define the cohomologies of $\Z^n$-graded Loday, Loday infinity,
and $2{p}$-ary graded Loday algebras. This leads to a graded Lie
``stem'' bracket, in which are encrypted, in addition to the
preceding cohomologies, the coboundary operators of graded Lie,
graded Poisson, graded Jacobi, Lie infinity, and $2{p}$-ary graded
Lie algebras \cite{MichVino}. Our approach is self-contained,
results are explicit, and they suggest that the operadic theory is
not confined to the quadratic case.
\medskip

The paper is organized as follows.\medskip

In Section 2, we study the properties of cohomology algebras,
which are implemented by square-zero elements in a graded Lie
algebra, with respect to formal deformations of these elements.
Our upshots extend similar properties for the adjoint Hochschild
(resp. Chevalley-Eilenberg, Leibniz) cohomology and deformations
of associative (resp. Lie, Loday) structures, which were proved in
\cite{Ger} (resp. \cite{NR}, \cite{Bal1}), and recovered in
\cite{Bal2}.
\medskip

Section 3 contains the definition of a graded dual Leibniz
coalgebra structure $\zD$ on the tensor algebra $T(W)$ of a
$\Z^n$-graded vector space $W$. We provide explicit formul{\ae}
for the reconstruction of coderivations and cohomomorphisms from
their corestriction maps.\medskip

In Section 4, we transfer the $\Z^n$-graded Lie bracket of
coderivations of the mentioned tensor coalgebra $T(W)$ of the
desuspension space $W:=\da V$ of an underlying $\Z^n$-graded
vector space $V$, and get a $\Z^{n+1}$-graded (resp.
$\Z^n$-graded) Lie bracket on the $\Z^{n+1}$-graded vector space
of weighted multilinear maps on $V$ (resp. on the $\Z^n$-graded
vector space of sequences of shifted weighted multilinear maps on
$V$). We determine the explicit form of this pullback ``stem''
bracket and show that its $\Z^{n+1}$-graded version coincides in
the case of a nongraded underlying space $V$ (resp. of graded
skew-symmetric multilinear mappings on $V$) with Rotkiewicz's
bracket \cite{Ro} pertaining to left Loday structures [and
corresponds to Balavoine's bracket \cite{Bal2} concerning right
Loday structures] (resp. with the graded Nijenhuis-Richardson
bracket \cite{LMS}).
\medskip

Codifferentials of our dual Leibniz coalgebra $T(W)$ are
characterized in Section 5. We prove that $\Z^n$-graded Loday
structures on $V$ can be viewed as (resp. we define strongly
homotopy Loday structures on $V$ as) degree
$e_1:=(1,0,\ldots,0)\in\Z^n$ quadratic (resp. odd degree)
codifferentials of $(T(\da V),\zD)$. Loday infinity structures and
Loday infinity morphisms are described in terms of sequences of
weighted multilinear maps that satisfy explicitly depicted
sequences of constraints: our $\op{Lod}_{\infty}$ algebras are
really differential graded Loday algebras up to homotopy and the
$\op{Lod}_{\infty}$ category contains the $\op{L}_{\infty}$
category as a subcategory.\medskip

Loday infinity (quasi)-isomorphisms are investigated. In Section
6, we prove a minimal model theorem for strongly homotopy Loday
algebras, and deduce that any Loday infinity quasi-isomorphism has
a quasi-inverse -- a theorem whose Lie infinity counterpart plays
a key-role in Deformation Quantization.\medskip

In Section 7, we deal with graded and strongly homotopy
cohomologies. $\Z^n$-graded Loday [resp. strongly homotopy Loday]
structures are canonical for the $\Z^{n+1}$-graded [resp.
$\Z^{n}$-graded] stem bracket, so that we obtain a natural
cohomology theory and an explicit coboundary operator. In the
nongraded (resp. the antisymmetric) [resp. the Lie infinity] case,
our $\Z^n$-graded Loday [resp. Loday infinity] cohomology operator
coincides with the Loday (resp. graded Chevalley-Eilenberg) [resp.
Lie infinity] differential given in \cite{DT} and \cite{Bal2}
(resp. in \cite{LMS}) [resp. in \cite{Pen} and \cite{FP}].

Further, graded Poisson and Jacobi cohomologies were defined
purely algebraically by Grabowski and Marmo in \cite{GM2}. The
authors prove existence and uniqueness of a $\Z^{n+1}$-graded
Jacobi (resp. Poisson) bracket on the algebra of antisymmetric
graded first order polydifferential operators (resp. of graded
polyderivations). We compute this ``Grabowski-Marmo'' bracket
explicitly and explain how the corresponding cohomologies are
induced by our stem bracket.

Finally, essentially two $p$-ary extensions of the Jacobi identity
were investigated during the last decades. The first, see e.g.
\cite{Fil}, leads to the Nambu-Lie structure, see \cite{Nam}, the
second, see \cite{MichVino}, \cite{VinoVino}, \cite{VinoVino2},
will in this text be referred to as $p$-ary Lie structure. We
define analogously $p$-ary ($p$ even) $\Z^n$-graded Loday
structures and their cohomology. These graded $p$-ary Loday
algebras are special strongly homotopy Loday algebras, so that we
have to prove that the two stem bracket induced cohomologies
coincide.

\section{Canonical elements of graded Lie algebras}\label{sectionCohomologyofagradedcanonicalelement}

Unless otherwise stated, all vector spaces that we consider in
this text are spaces over a field $\K$ of characteristic $0$, and
all graded vector spaces are $\Z^n$-graded, $n\in\N^*$. The
$\Z^n$--degree $\op{deg}(v)$ of a vector $v$ or the $\Z^n$--weight
$\op{deg}(f)$ of a graded linear map $f$ are often denoted by the
same symbol $v$ or $f$. If $v,f\in\Z^n$ are two such degrees, we
set $\la v,f\ra=\sum_iv_if_i.$  A homogeneous vector or graded
linear map $w$ is termed \emph{odd}, if $\la w,w\ra\in\Z$ is an
odd number. Eventually, again except differently stipulated, all
graded brackets are of weight $0$.

\begin{defi} We call \emph{canonical element} of a graded Lie algebra
(GLA) $(V,\{-,-\})$, any odd element $\zp\in V$ that verifies
$\{\pi,\pi\}=0$.
\end{defi}

It follows from the Jacobi identity that a GLA $(V,\{-,-\})$
endowed with a canonical element $\pi$ is a differential graded
Lie algebra (DGLA) $(V,\{-,-\},\p_{\pi})$ for
$\p_{\pi}:=\{\pi,-\}.$ Of course the cohomology space of any DGLA
is canonically a GLA. We denote the cohomology GLA of the
preceding DGLA by $H_{\pi}(V)$.\medskip

Let us mention two basic examples. If $V$ is a $\Z^n$-graded
vector space, we set $M(V) = \bigoplus_{(A,a)\in\Z^n\times\Z}
\linebreak M^{(A,a)}(V)$, where $M^{(A,a)}(V)=0$ for all $a\le
-2$, $M^{(A,-1)}(V)=V^A$, and where for each $a\ge 0$,
$M^{(A,a)}(V)$ is the space of all $(a+1)$-multilinear maps
$A:V^{\times(a+1)}\raa V$ that have weight $A$. Further, we denote
by $A(V)$ the graded vector subspace of $M(V)$, which is made up
by the $\Z^n$-graded skew-symmetric maps.

\begin{ex} The pair $(M(V),[-,-]^{\op{G}})$ (resp. $(A(V),[-,-]^{\op{NR}})$), where
$[-,-]^{\op{G}}$ (resp. $[-,-]^{\op{NR}}$) denotes the graded
Gerstenhaber (resp. Nijenhuis-Richardson) bracket, is a
$\Z^{n+1}$-graded Lie algebra. Its canonical elements of degree
$(0,1)\in\Z^n\times\Z$ are the associative graded (resp. graded
Lie) algebra structures on $V$. If $\zp$ is such an element,
$[\zp,\zp]^{\op{G}}=0$ (resp. $[\zp,\zp]^{\op{NR}}=0$) exactly
means that $\zp$ is associative (resp. verifies the graded Jacobi
identity). The cohomology $H_{\zp}(M(V))$ (resp. $H_{\zp}(A(V))$)
coincides with the adjoint Hochschild (resp. Chevalley-Eilenberg)
cohomology of $(V,\zp)$.\end{ex}

For details pertaining to these examples, we refer the reader to
\cite{LMS}.\medskip

Next we show that cohomology algebras $H_{\zp}(V)$ that are
implemented by canonical elements $\zp$ are good tools to study
formal deformations of $\pi$.\medskip

We set $V[[\nu]]=\bigoplus_{\za\in\Z^n}V^{\za}[[\nu]],$ where
$V^{\za}[[\nu]]$ is the space of formal power series in a formal
parameter $\nu$ with coefficients in the degree $\za$ term
$V^{\za}$ of $V$. A formal power series
$\label{FormalDeform}\pi_{\nu}=\sum_{i=0}^{\infty}\nu^i\pi_{i}\in
V^{\mathrm{deg}(\pi)}[[\zn]]$ with first term $\zp_0=\zp$ is a
\emph{formal deformation} of $\pi$, if it squares to zero w.r.t.
the natural extension of the bracket $\{-,-\}$ to a bilinear map
of the space $V[[\nu]]$, i.e. if
$\{\pi_{\nu},\pi_{\nu}\}=\sum_{p=0}^{\infty}\nu^p\sum_{i+j=p}\{\pi_i,\pi_j\}=0.$
A \emph{formal deformation of order $q$} is a formal series
$\zp_{\zn}$ that is truncated at order $q$ in $\nu$ and satisfies
the condition \be
\label{deforcondition}\sum_{i+j=p}\{\pi_i,\pi_j\}=0,\quad 1\leq
p\leq q.\ee We refer to formal deformations of order $1$ as
\emph{infinitesimal deformations}.

\begin{prop}
The cohomology space $H_{\pi}^{2\,\mathrm{deg}(\pi)}(V)$ contains
the obstructions to extension of formal deformations of order at
least $1$ to higher order deformations.
\end{prop}

\begin{proof} Assume that $\zp$ admits an order $q$, $q\geq 1$,
deformation $\pi_{\nu}$ and define $E_p$ in such a way that
Condition (\ref{deforcondition}) reads $E_{p}=-2\p_{\pi}(\pi_p),$
$1\leq p\leq q.$ In view of (\ref{deforcondition}), we then have
$$\p_{\pi}(E_{q+1})=-\sum_{\substack{k+l+j=q+1\\k,l,j\neq0}}\{\{\pi_k,\pi_l\},
\pi_j\},$$ which vanishes due to Jacobi's identity. In order to
extend deformation $\pi_{\nu}$ to order $q+1$, see
(\ref{deforcondition}), cocycle $E_{q+1}\in V^{2\op{deg}(\zp)}$
must be a coboundary.\end{proof}

\begin{lem}\label{Lemma:ExponMap} Let $\pi$ be a canonical element of a GLA
$(V,\{-,-\})$ and consider a series
$\chi_{\nu}=\sum_{i=1}^{\infty}\nu^i\chi_{i}\in V^0[[\zn]].$ If
$\pi_{\nu}$ is a formal deformation of $\pi$, then
$\op{exp}(\op{ad}\chi_{\nu})\;\pi_{\nu}$ is a formal deformation
of $\pi$ as well.
\end{lem}

\begin{proof} Let us mention that $\op{exp}$ denotes the
exponential series and stress that
\be(\op{ad}\chi_{\nu})^k\pi_{\nu}=\{\chi_{\nu}\{\chi_{\nu}\ldots\{\chi_{\nu},\pi_{\nu}\underbrace{\}\ldots\}\}}_k
=\sum_{p=0}^{\infty}\zn^p\sum_{i_1+\ldots+i_k+j=p}
\{\chi_{i_1}\{\chi_{i_2}\ldots\{\chi_{i_k},\zp_j\}\ldots\}\}.\label{ExpoSeries}\ee
It follows that the coefficient of $\zn ^p$ in the exponential
series over $k$ is made up by a finite number of terms in
$V^{\op{deg}(\zp)}$; indeed, if $k\ge p+1$, at least one of the
$\chi_{i_{\E}}$ vanishes. Moreover, the coefficient of $\zn^0$
contains only the term $k=0$, and thus
$(\op{exp}(\op{ad}\chi_{\nu})\;\pi_{\nu})_0=\zp$. Eventually, as
$\op{ad}\chi_{\nu}$ is a derivation of the bracket $\{-,-\}$, we
have
$$(\op{ad}\chi_{\nu})^{k}\{\pi_{\nu},\pi_{\nu}\}=\sum_{\substack{r+s=k\\r,s\geq0}}C_k^r\{(\op{ad}\chi_{\nu})^{r}\pi_{\nu},(\op{ad}\chi_{\nu})^{s}\pi_{\nu}\},$$
where $C_k^r$ is the binomial coefficient. Hence,
$\op{exp}(\op{ad}\chi_{\nu})\{\pi_{\nu},\pi_{\nu}\}=\{\op{exp}(\op{ad}\chi_{\nu})\pi_{\nu},\op{exp}(\op{ad}\chi_{\nu})\pi_{\nu}\},$
which completes the proof of the lemma.\end{proof}

Recall that two formal deformations $\pi_{\nu}$ and $\pi_{\nu}'$
of $\pi$ are said to be \emph{equivalent} (resp. \emph{equivalent
up to order $q$, $q\ge 1$}), if there is a series $\chi_{\zn}$ of
the type specified in Lemma \ref{Lemma:ExponMap}, such that
$\op{exp}(\op{ad}\chi_{\nu})\;\pi_{\nu}=\pi_{\nu}'$ (resp.
$\op{exp}(\op{ad}\chi_{\nu})\;\pi_{\nu}=\pi_{\nu}'+{\cal
O}(\zn^{q+1})).$ A deformation $\pi_{\nu}$ of $\pi$ is called
\emph{trivial} $($resp. \emph{trivial up to order} $q$,
$q\geq1$$)$, if $\pi_{\nu}$ is equivalent to $\pi$ $($resp.
equivalent to $\pi$ up to order $q$$)$.

\begin{prop} If $H_{\pi}^{\mathrm{deg}(\pi)}(V)=0$, any formal deformation of $\pi$ is trivial.
\end{prop}

\begin{proof} Let $\pi_{\nu}:=\pi+\sum_{i=1}^{\infty}\nu^i\pi_{i}$ be a
formal deformation of $\pi$. We first prove that $\pi_{\nu}$ is
trivial up to order 1, then we proceed by induction. Condition
(\ref{deforcondition}) and the assumption imply that there exists
$\chi_{1}\in V^0$, such that $\pi_1=\p_{\pi}(\chi_{1})$. When
setting $\chi_{\nu}^{(1)}=\nu\chi_{1}$, we get
$\op{exp}(\op{ad}\chi_{\nu}^{(1)})\pi_{\nu}=\pi+{\cal
O}(\nu^{2}).$ Suppose now that $\pi_{\nu}$ is trivial up to order
$q$ ($q\geq 1$), or, equivalently, that there is a series
$\chi_{\nu}^{(q)}$, such that
$\pi'_{\nu}:=\op{exp}(\op{ad}\chi_{\nu}^{(q)})\pi_{\nu}=\pi+
\nu^{q+1}\pi_{q+1}'+{\cal O}(\nu^{q+2}).$ As above, since
$\pi'_{\nu}$ is a deformation of $\pi$, there exists
$\chi_{q+1}'\in V^{0}$ that verifies
$\pi_{q+1}'=\p_{\pi}(\chi_{q+1}')$. Set now
$\chi_{\nu}^{(q+1)}:=\chi_{\nu}^{(q)}+\nu^{q+1}\chi_{q+1}'$ and
$\pi''_{\nu} :=\op{exp}(\op{ad}\chi_{\nu}^{(q+1)})\pi_{\nu}.$ It
follows from Equation (\ref{ExpoSeries}) that
$\pi''_{\nu}-\pi'_{\nu}=-\nu^{q+1}\p_{\pi}(\chi_{q+1}')+{\cal
O}(\nu^{q+2}).$ Hence, $\pi''_{\nu}=\pi+{\cal O}(\nu^{q+2})$,
which completes the proof. \end{proof}

Concerning infinitesimal deformations, it is easily seen from the
above explanations that

\begin{prop} Infinitesimal deformations of a canonical element
$\zp$ of a GLA $(V,\{-,-\})$ are classified up to first order
equivalence by $H^{\op{deg}(\zp)}_{\zp}(V)$.
\end{prop}

\section{Graded dual Leibniz tensor coalgebra}

Let us briefly recall some well-known facts. A \emph{graded
coalgebra} $(C,\Delta)$ is a graded vector space
$C=\bigoplus_{\alpha\in \Z^n}C^{\za}$ together with a
\emph{coproduct} $\Delta$, i.e. a linear map $\Delta:C\raa
C\otimes C$ that verifies
$\Delta(C^{\za})\subset\bigoplus_{\beta+\gamma=\za}
C^{\beta}\otimes C^{\gamma}.$  A {\em cohomomorphism} from
$(C,\zD)$ to a graded coalgebra $(C',\zD')$ is a weight $0$ linear
map ${\cal F}:C\to C'$, such that $\label{ComorphCond}
\Delta'{\cal F}=({\cal F}\otimes{\cal F})\Delta.$ In this text,
the tensor product of linear maps is defined by $(f\otimes
g)(v_1\otimes v_2)=(-1)^{\la g,v_1\ra}f(v_1)\otimes g(v_2),$ with
self-explaining notations. Further, a homogeneous
\emph{coderivation} of $(C,\Delta)$ is a linear map $Q:C\to C$ of
weight $\mathrm{deg}(Q)$ that satisfies the \emph{co-Leibniz
identity} $\label{coderivatonide} \Delta
Q=(Q\otimes\op{id}+\op{id}\otimes Q)\Delta,$ where $\op{id}$ is
the identity map of $C$. Weight $\za$ coderivations form a vector
space $\op{CoDer}^{\za}(C)$, and the space
$\op{CoDer}(C)=\bigoplus_{\za\in\Z^n}\op{CoDer}^{\za}(C)$ of all
coderivations carries a natural $\Z^n$-graded Lie algebra
structure provided by the graded commutator bracket.\medskip

To any $\Z^n$-graded vector space $V$, we associate the (reduced)
associative tensor algebra $T(V)=
\bigoplus_{p=1}^{\infty}V^{\otimes p}$ (the full tensor algebra
includes the term $V^{\otimes0}={\small \K}$ as well), which
carries two natural gradings, the $\Z$-gradation
$T(V)=\bigoplus_{p=1}^{\infty}T^pV$, $T^pV:=V^{\otimes p},$ and
the $\Z^n$-gradation $T(V)=\bigoplus_{\za\in\Z^n}T(V)^{\za}$,
$T(V)^{\za}:=\bigoplus_{p=1}^{\infty}(T^p{V})^{\za},$
$(T^p{V})^{\za}=\bigoplus_{\zb_1+\ldots+\zb_p=\za}V^{\zb_1}\ot\ldots\ot
V^{\zb_p}.$ In the following, unless differently stated, we view
$T(V)$ as $\Z^n$-graded vector space.\medskip

In order to define and study a coproduct on $T(V)$ that is adapted
to Loday structures, we need additional notations. For any
$p$--tuple $N^{(p)}:=(1,\ldots,p),$ $p\in\N^*$, an
\emph{unshuffle} $I=(i_1,\ldots,i_k)$, $1\leq k\leq p$, of
$N^{(p)}$ is a naturally ordered subset of $N^{(p)}$. The length
of $I$ is denoted by $|I|$. If $I$ and $J$ are two nonintersecting
unshuffles, we set
$(I;J)=(i_1,\ldots,i_{|I|};j_1,\ldots,j_{|J|})$, and $I\cup J$ is
the unique unshuffle that coincides with $(I;J)$ as a set.
Similarly,
$V_{(I;J)}=(v_{i_1},\ldots,v_{i_{|I|}};v_{j_1},\ldots,v_{j_{|J|}})$,
$v_{\E}\in V^{v_{\E}}$. The sign $(-1)^{(I;J)}$ is the signature
of the permutation $(I;J)\to I\cup J$, whereas $\ze_V(I;J)$ is the
Koszul sign implemented by $V_{(I;J)}\to V_{I\cup J}.$\medskip

\begin{prop}
Let $V$ be a $\Z^n-$graded vector space. The coproduct
$\Delta:T(V)\to T(V)\bigotimes T(V),$ defined by
\begin{equation}\label{coprodtensoralgebra}\Delta(v_1\otimes\ldots\otimes
v_p)=\sum_{\substack{I\cup
J=N^{(p-1)}\\I\neq\emp}}\varepsilon_V(I;J)\;V_I\bigotimes
V_J\otimes v_{p} \quad (v_{\E}\in V^{v_{\E}}, p\geq 1),
\end{equation} provides a graded dual Leibniz coalgebra structure on $T(V)$, i.e. a graded coalgebra structure that verifies
\begin{equation}\label{commultiplicationidentity}
    (\op{id}\bigotimes\Delta)\Delta=(\Delta\bigotimes\op{id})
    \Delta+(T\bigotimes\op{id})(\Delta\bigotimes\op{id})\Delta,
\end{equation}
where $T:T(V)\bigotimes T(V)\to T(V)\bigotimes T(V)$ is the
twisting map, which exchanges two elements of the $\Z^n$-graded
space $T(V)$ modulo the corresponding Koszul sign.
\end{prop}

% We would like to enlighten on why we have chosen an interest in this coalgebra structure, before we go
% ahead with laying out the Proof to this Proposition. Let $\pi=\{-,-\}$ be a Loday bracket. Seeing that the Jacobi
% identity is equivalent to \bea  \pi(\pi\bigotimes\op{id})=\pi(\op{id}\bigotimes\pi) +\pi(\op{id}\bigotimes\pi)(T\bigotimes\op{id}),\eea then,
% if $\delta$ denotes the dual of $\pi$, the restriction of  $\Delta$ to $V\ot V$ is exactly $\delta$.\\

\begin{proof} It suffices to compute the images of $v_1\otimes\ldots\otimes v_p$ by the three terms of
(\ref{commultiplicationidentity}). The first (resp. second; third)
one is equal to \be \sum_{\substack{I\cup K\cup
L=N^{(p-1)}\\I,K\neq\emp}}\varepsilon_V(I;K;L)V_I\bigotimes
V_K\bigotimes V_L\otimes v_p\ee (resp. to the same sum, but
confined to those unshuffles that verify $i_{|I|}<k_{|K|}$;
$k_{|K|}<i_{|I|}$, since $\varepsilon_V(K;I;L)(-1)^{\la
V_I,V_K\ra}=\varepsilon_V(I;K;L)$). Hence, the result. \end{proof}

%Subsequently, we shall say that two unshuffles $I,J$ verifies $I<
%J$ if $i_{|I|}<j_1$ that is to say $i_k<j_{\ell}$ for any
%$k\in\{1,\ldots,|I|\}$ and $\ell\in\{1,\ldots,|J|\}$

%\begin{rem}\label{remcoasscopro} Recall that the coassociative coalgebra structure on $T(V)$ (e.g. see \cite{St}) is given by the coproduct \be\label{coproductcoaass}\Delta_{\op{Ass}}(v_1\otimes\ldots\otimes v_p):=\sum_{i=1}^pv_1\ot\ldots\ot v_{i-1}\bigotimes v_{i}\ot\ldots\ot v_{p}, \quad\forall p\geq 1.\ee One easily ascertains that the part of (\ref{coprodtensoralgebra}) constituted by the sum of the terms of the forms $V_I\bigotimes V_J\ot v_p$ so that the unshuffles $I$ and $J$ satisfy $I<J$, coincides with $\Delta_{\op{Ass}}(v_1\otimes\ldots\otimes v_p)$.  \\
%\end{rem}

\begin{theo}\label{propositioncoderiavation}

The mapping \be\label{bijectioncod} \psi^V_Q:
\op{CoDer}^{Q}({T}(V))\ni Q\to
(Q_1,Q_2,\ldots)=:\sum_pQ_p\in\prod_{p\geq1}M^{(Q,p-1)}(V),\ee
which assigns to any (weight $Q$) coderivation $Q$ its (weight
$Q$) \emph{corestriction maps} $Q_p:T^pV\hookrightarrow
T(V)\xrightarrow{Q} T(V)\xrightarrow{\op{pr}}V,$ where $\op{pr}$
denotes the canonical projection, is a vector space isomorphism,
the inverse of which associates to any sequence $(Q_1,Q_2,\ldots)$
the coderivation $Q$ that is defined by \bean
\label{coderivationformula}
  &&Q(v_1\otimes\ldots\otimes v_{p})=\\ \nonumber &&\;\;\;\;\;\;
  \sum_{\substack{I\cup J\cup K=N^{(p)}\\I,J<K}}\varepsilon_V(I;J)(-1)^{\la Q,V_I\ra}V_I\otimes
  Q_{|J|+1}(V_J\otimes v_{k_1})\otimes V_{K\backslash k_1},\eean
where $I<K$ means that $i_{|I|}<k_1$ and where $v_{\E}\in
V^{v_{\E}}$.
\end{theo}

\begin{rem}\label{RemIsom1} The isomorphisms $\psi^V_Q$, $Q\in\Z^n$, (for inverses, see Equation (\ref{coderivationformula}))
induce an isomorphism $\psi^V$ between $\op{CoDer}(T(V))$ and the
corresponding direct sum of direct products. Further, if we denote
by $\op{CoDer}^Q_{p}(T(V))$, $Q\in\Z^n,p\in\N^*$, the image by
$(\psi^V_Q)^{-1}$ of $M^{(Q,p-1)}(V)$, isomorphism $\psi^V_Q$
restricts to an isomorphism $\psi^V_{(Q,p)}$ between these spaces.
If no confusion arises, we write $\psi$ instead of $\psi^V$,
$\psi^V_Q$, or $\psi^V_{(Q,p)}$.
\end{rem}

\begin{proof} It suffices to show that Equation
(\ref{coderivationformula}) defines a coderivation and that the
thus defined map $\Psi$ is the inverse of $\psi$. One easily sees
that $\psi\circ\Psi=\op{id}$. The second condition
$\Psi\circ\psi=\op{id}$ means that $\psi$ is injective. As
$\Delta$ vanishes if and only if its argument belongs to $V$, an
induction on $p$ allows deducing from the coderivation condition
and the definition of $\zD$ that $Q$ vanishes on $T^pV$, if its
corestriction maps vanish. As for the coderivation condition,
although it is quite easily checked for $p\le 3$, the general
proof is rather technical: it will not be reproduced here, but can
be found in \cite{APAXV}.\end{proof}

Like coderivations, cohomomorphisms from $(T(V),\zD)$ to
$(T(V'),\zD)$ are characterized by their corestriction maps.

\begin{rem}\label{RemCohomo} Let $V$ and $V'$ be two $\Z^n$-graded vector spaces. A coalgebra cohomomorphism ${\cal F}: (T(V),\Delta)\lraa
(T(V'),\Delta)$ is uniquely defined by its (weight 0)
corestriction maps $\cf_p:T^pV\to V'$, $p\geq1,$ via the equation
\be\label{Formu:CohoMorphism} \cf(v_1\otimes\ldots\otimes
v_{p})=\sum_{s=1}^{p}\sum_{\substack{I^1\cup\ldots\cup
I^s=N^{(p)}\\I^1,\ldots,
I^s\neq\emptyset\\i^{1}_{|I^1|}<\ldots<i^{s}_{|I^s|}}}\varepsilon_{V}(I^1;\ldots;I^s)\cf_{|I^1|}(V_{I^1})\ot\ldots\ot
\cf_{|I^s|}(V_{I^s}),\ee where $v_{\E}\in V^{v_{\E}}$. \end{rem}

Set now $e_1=(1,0,\ldots,0)\in\Z^n$ and consider the desuspension
operator $\da\,: V\raa\, \da V,$ where $\da V$ is the same space
as $V$ up to the shift $(\da V)^\za=V^{\za+e_1}$ of gradation. The
inverse map of $\da$ is denoted by $\uparrow$. The mapping
$\da^{\otimes p}:=\da\otimes\ldots\otimes\da,$ $p$ factors, i.e.
the mapping \be\label{signeruledesuspention} \da^{\otimes p}:
V^{\otimes p}\ni v_1\otimes\ldots\otimes v_p\to
(-1)^{\sum_{s=1}^p\la (p-s)e_1,v_s\ra}\downarrow
v_1\otimes\ldots\otimes\downarrow v_p\in(\downarrow V)^{\otimes
p},\ee is an isomorphism of weight $-p\,e_1$, whose inverse is
$\label{desuspensioninverse}(-1)^{\frac{p(p-1)}{2}}\uparrow^{\otimes
p}.$

\begin{rem}\label{RemIsom2} The isomorphisms \be
\zs^{\da V}_{(Q,p)}:M^{(Q,p-1)}(\da V)\ni Q_p\to \zp_p:= \ua\circ
Q_p\circ\da^{\otimes p}\in
M^{(Q+(1-p)e_1,p-1)}(V),\label{EquIsom2}\ee $Q\in\Z^n,p\in\N^*$,
(their inverses are defined by $(-1)^{\frac{p(p-1)}{2}}\da\circ
\zp_p\circ \ua^{\otimes p}$) generate isomorphisms $\zs^{\da V}_Q$
and $\zs^{\da V}$ between the corresponding direct products and
direct sums of direct products. If no confusion is possible, we
omit super- and subscripts and denote these isomorphisms simply by
$\zs$. Isomorphisms (\ref{EquIsom2}) extend of course to
multilinear maps on $\da V$ valued in $\da V'$.
\end{rem}

\begin{rem}\label{RemSequenceView} Theorem \ref{propositioncoderiavation} and Remarks \ref{RemIsom1}-\ref{RemIsom2} show that weight $Q$
coderivations $Q:(T(\da V),\Delta)\raa(T(\da V),\Delta)$ can be
viewed as sequences $\zp=(\zp_1,\zp_2,\ldots)=\sum_p\zp_p$ of
weight $Q+(1-p)e_1$ multilinear maps $\zp_p:V^{\times p}\to V$,
and that cohomomorphisms $\cf:(T(\da V),\Delta)\raa(T(\da
V'),\Delta)$ (which by definition have weight $0$) can be seen as
sequences $f=(f_1,f_2,\ldots)=\sum_pf_p$ of weight $(1-p)e_1$
multilinear maps $f_p:V^{\times p}\to V'$.\end{rem}

\section{Stem bracket}\label{sectionGeneralizedgradedLiebracket}

When combining the isomorphisms $\zs^{-1}$ and $\psi^{-1}$, we
get, for $A\in\Z^n, a\in\N$, a vector space isomorphism
\be\phi_{(A,a)}: M^{(A,a)}(V)\ni A\to
Q^A=(0,\ldots,0,Q_{a+1}^A,0,\ldots)\in\op{CoDer}_{a+1}^{A+ae_1}(T(\da
V)),\ee where $Q_{a+1}^A=(-1)^{\frac{a(a+1)}{2}}\da\circ
A\circ\uparrow^{\ot(a+1)}$ and where $Q^A$ is the coderivation
that is obtained from $Q^A_{a+1}$ via extension equation
(\ref{coderivationformula}).

\begin{theo}\label{proplodaybrackets} The $\Z^{n+1}$-graded vector space $M_r(V)=M(V)\ominus V$ is a $\Z^{n+1}$-graded Lie
algebra, when endowed with the bracket \be\label{16}[A,B]^{\ot}
:=(-1)^{1+\la
ae_1,be_1+B\ra}\phi_{{(A+B,a+b)}}^{-1}([\phi_{(A,a)}(A),\phi_{(B,b)}(B)]),\ee
$A\in M^{(A,a)}(V)$, $B\in M^{(B,b)}(V)$, where $[-,-]$ denotes
the $\Z^n$-graded Lie bracket of the space of coderivations of
$(T(\da V),\zD)$.
\end{theo}

\begin{proof} It follows from Equation
(\ref{coderivationformula}) that the $p$-th corestriction map
$[\phi_{(A,a)}(A),\phi_{(B,b)}(B)]_p$ vanishes if $p\neq a+b+1$,
so that the bracket $[\phi_{(A,a)}(A),\phi_{(B,b)}(B)]$ is really
a coderivation in $\op{CoDer}^{A+B+(a+b)e_1}_{a+b+1}(T(\da V))$.
The sign $(-1)^{1+\la ae_1,be_1+B\ra}$ ensures that the
$\Z^n$-graded Lie bracket of coderivations induces a
$\Z^{n+1}$-graded Lie bracket $[-,-]^{\ot}$.\end{proof}

As the map $\zf=\psi^{-1}\circ\zs^{-1}$ is also a (weight 0)
$\Z^n$-graded vector space isomorphism \be \zf:
C(V):=\bigoplus_{Q\in\Z^n}\prod_{p\ge
1}M^{(Q+(1-p)e_1,p-1)}(V)\to\op{CoDer}(T(\da
V))=\bigoplus_{Q\in\Z^n}\op{CoDer}^Q(T(\da V)),\ee the next
proposition is obvious.

\begin{prop}\label{Zn+1Stem} The $\Z^n$-graded vector space $C(V)$ is a $\Z^n$-graded
Lie algebra for the bracket
\be[\zp,\zr]^{\bar{\otimes}}=\zf^{-1}[\zf\zp,\zf\zr]=\sum_{q\ge
1}\sum_{s+t=q+1}(-1)^{1+(s-1)\la e_1,\zr
\ra}[\zp_s,\zr_t]^{\otimes},\ee where $\zp=\sum_s\zp_s\in
C^{\zp}(V)$ and $\zr=\sum_t\zr_t\in C^{\zr}(V)$ are two
homogeneous $C(V)$-elements of $\Z^n$-degree $\zp$ and $\zr$
respectively.
\end{prop}

\begin{rem} In the following, we refer to $[-,-]^{\bar{\otimes}}$ (resp.
$[-,-]^{\otimes}$) as the $\Z^n$-graded (resp. $\Z^{n+1}$-graded)
stem bracket.\end{rem}

\begin{theo}\label{ExplLodaBrac} The $\Z^{n+1}$-graded stem bracket
on $M_r(V)$ explicitly reads \be\label{bracketloday}
[A,B]^{\ot}=j_AB-(-1)^{\la (A,a),(B,b)\ra}j_BA,\ee where
\bean\label{MLodaycomposition} &&(j_BA)(v_1\ot\ldots\ot
v_{a+b+1})=(-1)^{\la A,B\ra}\sum_{\substack{I\cup J\cup
K=N^{(a+b+1)} \nonumber\\I,J< K,\,|J|=b}}(-1)^{\la
B,V_I\ra+b|I|}\\&&\;\;\;\;\;\;\;\;\;\;\;\;
(-1)^{(I;J)}\varepsilon_{V}(I;J)A(V_I\otimes B(V_J\otimes
v_{k_1})\otimes V_{K\backslash k_1}),\eean for any $A\in
M^{(A,a)}(V)$, $B\in M^{(B,b)}(V)$, $a,b\ge 0$, and $v_{\E}\in
V^{v_{\E}}$.
\end{theo}

\begin{proof} It follows from Equation
(\ref{coderivationformula}) that the $(a+b+1)$-th corestriction of
$[\phi_{(A,a)}(A),\phi_{(B,b)}(B)]=[Q^A,Q^B]$ is obtained by just
restricting the involved composite maps to $(\da
V)^{\otimes(a+b+1)}.$ The description of the isomorphisms
$\zf_{(A,a)}^{-1}$ then shows that \bean\label{LodaBrac}
[A,B]^{\ot}&=&\nonumber -(-1)^{\la
ae_1,be_1+B\ra}\left(\uparrow\circ Q^A\circ
Q^B\circ\da^{\otimes(a+b+1)}\right.\\&&\left.-(-1)^{\la
A+ae_1,B+be_1\ra}\uparrow\circ Q^B\circ Q^A
\circ\da^{\otimes(a+b+1)}\right).\eean

If $V_{N^{(a+b+1)}}=v_1\ot\ldots\ot v_{a+b+1},$ with $v_{\E}\in
V^{v_{\E}}$, we get \bea \left(\uparrow\circ Q^A\circ
Q^B\circ\da^{\otimes(a+b+1)}\right)(V_{N^{(a+b+1)}})&=&(-1)^{\zb_1}\left(\uparrow\circ
Q^A\circ Q^B\right)(\da V_{N^{(a+b+1)}}),\\\zb_1&=&\la
e_1,\sum_{s\ge 1}(a+b+1-s)v_s\ra,\eea where $\da
V_{N^{(a+b+1)}}=\da v_1\otimes\ldots\ot\da v_{a+b+1}$. Formula
(\ref{coderivationformula}) yields \bean Q^B(\da
V_{N^{(a+b+1)}})&=&\sum_{\substack{I\cup J\cup
K=N^{(a+b+1)}\\I,J<K,|J|=b}}(-1)^{\zb_2}\varepsilon_{\da
V}(I;J)\da V_I\otimes Q^{B}_{b+1}(\da V_J\otimes \da
v_{k_1})\otimes \da V_{K\backslash k_1},\nonumber\\\zb_2&=&\la
B+be_1,V_I+|I|e_1\ra.\eean Moreover, \bea Q^{B}_{b+1}(\da
V_J\otimes \da v_{k_1})= (-1)^{\zb_3}\da B(V_J\otimes
v_{k_1}),\zb_3=\la e_1,\sum_{s\ge 1}(b+1-s)v_{j_s}\ra,\eea as the
sign $(-1)^{\frac{b(b+1)}{2}}$ inside $Q^B_{b+1}$ and the sign due
to the shift of the $\Z^n$-gradation cancel each other out. When
noticing that $Q^A$ evaluated on an element of $(\da
V)^{\otimes(a+1)}$ is nothing but $Q_{a+1}^A$, we obtain

\bean && \label{17}\left(\uparrow\circ Q^A\circ
Q^B\circ\da^{\otimes(a+b+1)}\right)(V_{N^{(a+b+1)}})=\\&&\;\;\;\;\;\;\;\;\;\;\;\;\;\;\sum_{\substack{I\cup
J\cup K=N^{(a+b+1)}\\I,J<K,\;|J|=b}}\nonumber
(-1)^{\ell}\varepsilon_{\da V}(I;J) A(V_I\otimes B(V_J\otimes
v_{k_1})\otimes V_{K\backslash k_1}),\eean with
$\ell=\zb_1+\zb_2+\zb_3+\zb_4+\zb_5+\zb_6,$ where $$\zb_4=\la
e_1,\sum_{s\ge 1}(a+1-s)v_{i_s}\ra,\;\zb_5=(a-\vert I\vert)\la
e_1, B+V_J+v_{k_1}\ra,\;\mbox{and}\,\;\zb_6=\la e_1,\sum_{s\ge
2}(a-\vert I\vert -s+1)v_{k_s}\ra$$ are generated by
$\uparrow^{\otimes(a+1)}$ and where, again, the sign inside
$Q^A_{a+1}$ and the sign due to the shift neutralize.

We will prove that \be\label{18}(-1)^{\ell}\varepsilon_{\da
V}(I;J)=(-1)^{\la B,ae_1\ra+\la
B,V_I\ra+b|I|}(-1)^{(I;J)}\varepsilon_{V}(I;J).\ee Observe first
that an appropriate regrouping of terms yields \bea
&&\ell=\E_1+\E_2:=\left(\la B,ae_1\ra+\la
B,V_I\ra+b|I|\right)+\left(\la e_1,\sum_{s=1}^{a+b+1}
(a+b+1-s)v_{s}\ra\right.\\\\&&+\la
e_1,\sum_{s=1}^{|I|}(a+b+1-s)v_{i_s}\ra+ \la
e_1,\sum_{s=|I|+1}^{|I|+|J|} (a+b+1-s)v_{j_{s-|I|}}\ra\\
&&\;\;\;\;\;\;\;\;\;\;\;\;\;\;\;\;\;\left.+ \la
e_1,\sum_{s=|I|+|J|+1}^{a+b+1}(a+b+1-s)v_{k_{s-|I|-|J|}}\ra\right),
\eea where, since in view of the conditions $I,J<K$ the
concatenation $(I;J)$ is a permutation of $1,\ldots,\vert
I\vert+\vert J\vert$, the term $\E_2$ reads (modulo even terms)
$$\E_2=\la e_1,\sum_{s=1}^{\vert I\vert}(s+i_s)v_{i_s}\ra+\la
e_1,\sum_{s={\vert I\vert}+1}^{\vert I\vert +\vert
J\vert}(s+j_{s-\vert I\vert})v_{j_{s-\vert I\vert}}\ra.$$ If
permutation $(I;J)$ is a transposition
$(I;J)=(1,\ldots,q-1,q+1,q,q+2,\ldots,\vert I\vert+\vert J\vert)$,
then $\E_2=\la e_1,v_q+v_{q+1}\ra.$ It is now easily checked that
for a transposition \be(-1)^{\E_2}\ze_{\da
V}(I;J)=(-1)^{(I;J)}\ze_V(I;J)\label{18bis}\ee and that for a
composition of transpositions all the factors of this last
equation are the products of the corresponding factors induced by
the involved transpositions. It follows that Equations
(\ref{18bis}) and (\ref{18}) hold true for any permutation
$(I;J)$.

Eventually Equation (\ref{17}) may be written \bean\nonumber && \left(\label{41}\uparrow\circ Q^A\circ Q^B\circ\da^{\otimes(a+b+1)}\right)(V_{N^{(a+b+1)}})=\\
&&\;\;\;\;\sum_{\substack{I\cup J\cup K=N^{(a+b+1)}\\I,J<
K,\;|J|=b}}\nonumber (-1)^{m}(-1)^{(I;J)}\varepsilon_{V}(I;J)
A(V_I\otimes B(V_J\otimes v_{k_1})\otimes V_{K\backslash
k_1}),\nonumber\eean where $m=\la B,ae_1\ra+\la B,V_I\ra+b|I|$.

If we define the insertion operator \benn j_BA=(-1)^{\la ae_1
+A,B\ra}\uparrow\circ Q^A\circ Q^B\circ\da^{\otimes(a+b+1)},\eenn
$A\in M^{(A,a)}(V)$, $B\in M^{(B,b)}(V)$, $a,b\ge 0$, we finally
get the announced result.\end{proof}

\begin{rem} \label{Rem:StemToNR}
It is easily checked that the restriction of the stem bracket
$[-,-]^{\ot}$ to the subspace $\op{A}_r(V)$ of $M_r(V)$, made up
by graded skew-symmetric multilinear maps on $V$, coincides with
the graded Nijenhuis-Richardson bracket $[-,-]^{\op{NR}}$ that was
introduced in \cite{LMS}.\end{rem}

\section{Graded and strongly homotopy Loday structures}

Let us recall that a \emph{codifferential} of a coalgebra is a
coderivation that squares to $0$.

\begin{prop}\label{procodifferentailcarac}
A homogenous odd weight coderivation $Q$ of the coalgebra
$(T(V),\zD)$ is a codifferential if and only if, for any $p\geq1$,
the following equation holds identically:

\be \label{codifferentialformula} \sum_{\substack{I\cup J\cup
K=N^{(p)}\\I,J<K}}\varepsilon_V(I;J)(-1)^{\la
Q,V_I\ra}Q_{|I|+|K|}(V_I\otimes Q_{|J|+1}(V_J\otimes
v_{k_1})\otimes V_{K\backslash k_1})=0.\ee
\end{prop}

\begin{proof} As $Q$ is an odd weight coderivation, $2Q^2=[Q,Q]$ and so
$Q^2$ is also a coderivation. Thus, according to Theorem
\ref{propositioncoderiavation}, the condition $Q^2=0$ is satisfied
if and only if the corestriction maps $Q^2_p$, $p\geq1$, of $Q^2$
vanish. It is easily seen that $Q^2_p(V_{N^{(p)}})$ is exactly the
LHS of (\ref{codifferentialformula}). \end{proof}

\begin{defi} A \emph{graded Loday algebra} $(V,\{-,-\})$ $($GLodA for
short$)$ is made up by a $\Z^n$-graded vector space $V$ and a
weight $0$ bilinear bracket $\{-,-\}$ that satisfies the graded
Jacobi identity \begin{equation}\{a, \{b, c\}\}=\{\{a,b\},c\} +
(-1)^{\la a,b\ra}\{b, \{a, c\}\},
\end{equation}
for any homogeneous $a,b,c\in V$.
\end{defi}

Nongraded Loday algebras were introduced by Loday in \cite{Lo1}
and are also known as {\em Leibniz algebras}. For further
information on graded Loday algebras, we refer the reader to
\cite{Kos} or \cite{Akm}.

\begin{theo}\label{lodaycodifferential}
Let $\op{Lod}(V)$ be the set of $\Z^n$-graded Loday structures on
a $\Z^n$-graded vector space $V$, and denote by
$\op{CoDiff}_{p}^Q(T(\da V))$, $Q\in\Z^n$, $p\in\N^*$, the set of
codifferentials $Q$ of $(T(\da V),\zD)$, which have weight $Q$ and
whose corestriction maps all vanish except $Q_{p}.$ Then,
\begin{enumerate}
  \item The restriction of $\phi_{(0,1)}$ to $\op{Lod}(V)$ is a
  bijection and
\begin{equation}\label{Lodaycodifferentialbijection}
\op{Lod}(V)\simeq \op{CoDiff}_{2}^{e_1}(T(\downarrow
V)),\end{equation}
  \item The graded Loday structures on $V$ are exactly the canonical elements of weight $(0,1)$ of the graded Lie algebra
  $(M_r(V),[-,-]^{\ot})$.
\end{enumerate}
\end{theo}

\begin{proof} Since $\phi_{(0,1)}$ is a bijection between $M^{(0,1)}(V)$
and $\op{CoDer}_{2}^{e_1}(T(\da V))$, it suffices, in order to
account for Point 1, to prove that
$\phi_{(0,1)}(\op{Lod}(V))=\op{CoDiff}_{2}^{e_1}(T(\downarrow
V)).$

If $\pi\in\op{Lod}(V)$, its image
$\phi_{(0,1)}(\pi)=Q^{\zp}=(0,Q^{\zp}_2,0,
\ldots)\in\op{CoDer}^{e_1}_2(T(\da V))$ is a codifferential, if
\bean\label{13} Q_{2}^{\pi}(Q_{2}^{\pi}(\da v_1\otimes \da
v_2)\otimes \da v_3)+(-1)^{\la e_1, \da v_1\ra}\nonumber
Q_{2}^{\pi}(\da v_1\otimes Q_{2}^{\pi}(\da v_2\otimes \da
v_3))\\\;\;\;\;\;\;\;\;\;\;\;\;\;\;\;\;\;\;\;\;+(-1)^{\la e_1+\da
v_1,\da v_2\ra}Q_{2}^{\pi}(\da v_2\otimes Q_{2}^{\pi}(\da
v_1\otimes \da v_3))=0,\eean see Proposition
\ref{procodifferentailcarac} and note that Condition
(\ref{codifferentialformula}) is trivial for $p\neq 3.$ When
remembering that $Q_{2}^{\pi}=-\da\circ\pi\circ
(\uparrow\otimes\uparrow)$, we easily check that Condition
(\ref{13}) reads \bea &&(-1)^{\la v_2,e_1\ra}\da
[\pi(\pi(v_1\otimes v_2)\otimes
v_3)\\&&\;\;\;\;\;\;\;\;\;\;\;\;\;\;\;\;\;\;\;\;\;\;\;\;-\pi(v_1\otimes\pi
(v_2\otimes v_3))+(-1)^{\la v_1,v_2\ra}\pi(v_2\otimes\pi
(v_1\otimes v_3))]=0.\eea As $\pi$ verifies the graded Jacobi
identity, the last requirement is fulfilled.

Conversely, if
$Q=(0,Q_2,0,\ldots)\in\op{CoDiff}_{2}^{e_1}(T(\downarrow V))$,
then $\pi:=\phi^{-1}_{(0,1)}(Q)$ is a graded Loday structure.

As regards Point 2, note that any graded Loday structure $\pi$ is
odd. Furthermore, \bea [\pi,\pi]^{\ot}= -(-1)^{\la
e_1,e_1\ra}\phi_{(0,2)}^{-1}([\phi_{(0,1)}(\pi),\phi_{(0,1)}(\pi)])=2\phi_{(0,2)}^{-1}(\phi_{(0,1)}^2(\pi)),\eea
so that the graded Loday structures are exactly the canonical
elements of weight $(0,1)$.\end{proof}

There are two ways to make an algebraic structure more flexible,
categorification and homotopyfication. In the sequel, we extend
graded Loday structures, see Equation
(\ref{Lodaycodifferentialbijection}), and investigate the homotopy
version of Loday algebras.

\begin{defi}\label{Defi:StronLodayStrucCoal}
A {\em strongly homotopy Loday algebra} or a {\em Loday infinity
algebra} is a $\Z^n$-graded vector space $V$ endowed with a
codifferential of odd weight of the tensor coalgebra $(T(\da
V),\Delta)$.\end{defi}

We denote by $\op{Lod}^Q_{\infty}(V)$, $Q\in\Z^n$, $\la Q,Q\ra$
odd, the set
$$\op{Lod}_{\infty}^Q(V)\simeq\op{CoDiff}^Q(T(\da V))$$ of weight $Q$ Loday infinity
($\op{Lod}_{\infty}^Q$ for short) structures on $V$, whereas
$\op{Lod}_{\infty}(V)$ denotes the set
$\op{Lod}_{\infty}^{e_1}(V)$ -- as most infinity structures
considered below have weight $e_1$. Since $Q$ is odd,
$Q\in\op{Lod}^Q_{\infty}(V)$ if and only if
$Q\in\op{CoDer}^Q(T(\da V))$ and $[Q,Q]=2Q^2=0$. Hence, in view of
Remark \ref{RemSequenceView} and Proposition \ref{Zn+1Stem}, the
sequence ``definition'' of Loday infinity algebras:

\begin{prop} A Loday infinity algebra is a $\Z^n$-graded vector
space $V$ together with a sequence of \emph{structure maps}
$$\zp=(\zp_1,\zp_2,\ldots)=\sum_p\zp_p\in C^{Q}(V)=\prod_{p\ge 1}M^{(Q+(1-p)e_1,p-1)}(V)$$
of odd degree $Q$, such that \be \sum_{s+t=p}(-1)^{1+(s-1)\la
e_1,Q \ra}[\zp_s,\zp_t]^{\otimes}=0, \forall p\ge
2.\label{LodaInftCond}\ee
\end{prop}

In the usual case of $\op{Lod}_{\infty}$ structures $\zp$ on $V$,
the first three conditions (\ref{LodaInftCond}) mean that
$(V,\pi_1)$ is a chain complex, that $\pi_1$ is a $\Z^n$-graded
derivation of the bilinear map $\pi_2$, and that $\pi_2$ is a
$\Z^n$-graded Loday structure modulo homotopy $\pi_3$.

\begin{ex}
If the structure maps of a $\op{Lod}_{\infty}$ algebra $(V,\pi)$
all vanish, except $\zp_1$ $($resp. except $\zp_2$, except $\pi_1$
and $\pi_2$$)$, $(V,\zp)$ is a chain complex $($resp. a
$\Z^n$-graded Loday algebra (GLodA), a differential graded Loday
algebra$ (DGLodA))$.
\end{ex}

\begin{ex} \label{Ex:DirectSumLodayInfty} Let $(V,\pi)$ and $(V^{\prime},
\pi^{\prime})$ be two $\op{Lod}_{\infty}$ algebras. Their direct
sum $(V\oplus V^{\prime},\pi\oplus\pi^{\prime})$ is a
$\op{Lod}_{\infty}$ algebra, where the structure maps are defined
by
$$(\pi\oplus\pi^{\prime})_{p}(v_1+v_{1}^{\prime},\ldots,v_p+
v_{p}^{\prime}):=\pi_p(v_1,\ldots,v_p)+\pi_p^{\prime}(v^{\prime}_1,\ldots,v^{\prime}_p),$$
for any $p\geq1$.
\end{ex}

It is obvious from the structure of the explicit form of the stem
bracket $[-,-]^{\otimes}$, see Theorem \ref{ExplLodaBrac}, that if
$\zp$ and $\zp'$ verify the $\op{Lod}_{\infty}$ conditions
(\ref{LodaInftCond}) in $V$ and $V'$ respectively, then $\zp\oplus
\zp'$ verifies the same conditions in $V\oplus V'$.\medskip

As a $\op{Lod}_{\infty}$ structure on $V$ is a weight $e_1$
codifferential $Q$ of the coalgebra $(T(\da V),\zD)$, it is
natural to give the following coalgebraic

\begin{defi}\label{Defi:MorphismLoday}
Let $(V,Q)$ and $(V',Q')$ be two $\op{Lod}_{\infty}$ algebras. A
{\em $\op{Lod}_{\infty}$ morphism} from $(V,Q)$ to $(V',Q')$ is a
coalgebra cohomomorphism ${\cal F}: (T(\da V),\Delta)\lraa (T(\da
V'),\Delta),$ such that $Q'\cf=\cf Q.$
\end{defi}

Since the composition of $\op{Lod}_{\infty}$ morphisms is again a
$\op{Lod}_{\infty}$ morphism, and as, for any $\op{Lod}_{\infty}$
algebra, the identity map is a morphism (say ${\cal I}$, with
corestriction maps ${\cal I}_1=\op{id}$ and ${\cal I}_p=0$, for
$p\geq2$), $\op{Lod}_\infty$ algebras and their morphisms form a
category, which we denote by $\mathbf{Lod}_\infty$.\medskip

Next we give the sequence ``definition'' of $\op{Lod}_{\infty}$
morphisms, see Remark \ref{RemSequenceView}.

\begin{prop}
Let $(V,\pi)$ and $(V',\pi')$ be two $\op{Lod}_{\infty}$ algebras.
A $\op{Lod}_{\infty}$ morphism $f:(V,\pi)\raa (V',\pi')$ is a
sequence of weight $(1-p)e_1$ multilinear maps $f_p:V^{\times
p}\to V'$, $p\ge 1$, which satisfy the condition \bean\label{75}
&&\sum_{s=1}^{p}\sum_{\substack{I^1\cup\ldots\cup
I^s=N^{(p)}\\I^1, \ldots,
I^s\neq\emptyset\\i^{1}_{|I^1|}<\ldots<i^{s}_{|I^s|}}}(-1)^{\omega}(-1)^{(I^1;\ldots;I^s)}\nonumber
\varepsilon_V(I^1;\ldots;I^s)\pi'_s(f_{|I^1|}(V_{I^1}), \ldots,
f_{|I^s|}(V_{I^s}))\\&&=\sum_{\substack{I\cup J\cup
K=N^{(p)}\\I,J<K}}(-1)^{\lambda}(-1)^{(J;I)}\varepsilon_{V}(I;J)f_{|I|+|K|}(V_I,
\pi_{|J|+1}(V_J,v_{k_1}), V_{K\backslash k_1}),\eean where
\be\label{1981}\begin{array}{l}\omega=\frac{s(s-1)}{2}+\sum_{1\leq
r\leq s}(s-r)|I^{r}|+\sum_{2\leq r\leq s}\la(|I^{r}|+1)e_1,
V_{|I^{1}|}+\ldots+V_{|I^{r-1}|}\ra\end{array}\ee and
\be\lambda=\la(1+|J|)e_1,V_I+(p+1)e_1\ra,\ee for all $p\ge 1$.
\end{prop}

\begin{proof} When using (\ref{coderivationformula}) and
(\ref{Formu:CohoMorphism}), we see that $Q'\cf$ and $\cf Q$ are
similarly composed of the respective corestriction maps, so that
they coincide, if their corestrictions do. This leads to a
reformulation of the intertwining condition $Q'\cf=\cf Q$ in terms
of $Q'_p,{\cal F}_p,Q_p$, $p\ge 1$. The translation of this
reformulation by means of the sequences $\zp',f,\zp$, and in
particular the above signs, are obtained by a direct computation
that is based upon similar arguments than those used in the proof
of Theorem \ref{ExplLodaBrac}. We refrain from providing the
details.\end{proof}

\begin{rem}
The first constraint $(\ref{75})$ states that $f_1$ is a chain map
between $(V,\pi_1)$ and $(V',\pi'_1)$, whereas the second means
that $f_2$ measures the deviation from $f_1$ being a
$(V,\zp_2)-(V',\zp'_2)$ homomorphism. If $(V,\pi)$ and $(V',\pi')$
are DGLodAs, map $f_1$ is a DGLodA morphism. It can be shown that
the category $\mathbf{L}_{\infty}$ of $\op{L}_{\infty}$ algebras
and morphisms is a subcategory of
$\mathbf{Lod}_{\infty}$.\end{rem}

\begin{defi}
A {\em $\op{Lod}_{\infty}$ quasi-isomorphism} from a
$\op{Lod}_{\infty}$ algebra $(V,\pi)$ to a $\op{Lod}_{\infty}$
algebra $(V',\pi')$ is a $\op{Lod}_{\infty}$ morphism
$f:(V,\pi)\raa (V',\pi')$, such that the chain map
$f_1:(V,\pi_1)\raa (V',\pi'_1)$ induces an isomorphism
$f_{1\;\sharp}:H(V,\pi_1)\raa H(V',\pi_1')$ in cohomology. In
particular, $f$ is called a {\em $\op{Lod}_{\infty}$ isomorphism},
if $f_1:V\raa V'$ is an isomorphism.
\end{defi}

If $\cf\simeq f$ and $\cg\simeq g$ are two composable
$\op{Lod}_{\infty}$ morphisms, we denote by $g\circ f$ the
sequence of multilinear maps that corresponds to the
$\op{Lod}_{\infty}$ morphism $\cg\cf$. Similarly, $\pi'\circ f$
and $f\circ\pi$ are the sequences that represent $Q'{\cal F}$ and
${\cal F}{Q}$. The $\op{Lod}_{\infty}$ morphism condition
(\ref{75}) then reads $\pi'\circ f=f\circ\pi.$ We use these and
analogous notations below.\medskip

The next upshots will be needed in the following.

\begin{prop}\label{InveLodInftIsom}\begin{enumerate} \item Any coalgebra cohomomorphism $f:(T(\da V),\zD)\to (T(\da
V'),\zD)$, which corresponds to a sequence $f=(f_1,f_2,\ldots)$,
whose first element $f_1$ is bijective, is invertible, i.e. there
is a coalgebra cohomomorphism $f^{-1}:(T(\da V'),\zD)\to (T(\da
V),\zD)$, such that $f\circ f^{-1}=\op{Id}$ and $f^{-1}\circ
f=\op{Id}$, where $\op{Id}$ is the unit cohomomorphism
$\op{Id}=(\op{id},0,0,\ldots).$

\item If $(V,\zp)$ denotes a $\op{Lod}_{\infty}$ algebra, any
sequence $f=(f_1,f_2,\ldots)$ of weight $(1-p)e_1$ multilinear
maps $f_p:V^{\times p}\to V$, whose first element $f_1$ is the
identity map of $V$, induces a new $\op{Lod}_{\infty}$ structure
$f\circ \zp\circ f^{-1}$ on $V$ and $f$ is a $\op{Lod}_{\infty}$
isomorphism between $(V,\zp)$ and $(V,f\circ \zp\circ f^{-1})$.

\item Any $\op{Lod}_{\infty}$ isomorphism $f:(V,\zp)\to (V',\zp')$
admits an inverse $f^{-1}$ that is a $\op{Lod}_{\infty}$
isomorphism as well.\end{enumerate}\label{InverMorph}\end{prop}

\begin{proof} 1. Let ${\cal F}:(T(\da V),\zD)\to (T(\da V'),\zD)$ be a
coalgebra cohomomorphism, whose first corestriction ${\cal
F}_1:\da V\to \da V'$ is bijective. If there is an inverse
cohomomorphism ${\cal G}:(T(\da V'),\zD)\to (T(\da V),\zD)$, it
follows from the condition ${\cal F G}={\cal I}$ and from Equation
(\ref{Formu:CohoMorphism}) that ${\cal G}_1={\cal F}_1^{-1}$ and
that, for any $p\ge 2$, \bea &&{\cal G}_p(\da v'_1,\ldots,\da
v'_p)\\&&=-\sum_{s=2}^p\sum_{\substack{I^1\cup\ldots\cup
I^s=N^{(p)}\\I^1,\ldots,
I^s\neq\emptyset\\i^{1}_{|I^1|}<\ldots<i^{s}_{|I^s|}}}\varepsilon_{\da
V'}(I^1;\ldots;I^s){\cal F}_1^{-1}{\cal F}_s\lp{\cal
G}_{|I^1|}(\da V'_{I^1})\ot\ldots\ot {\cal G}_{|I^s|}(\da
V'_{I^s})\rp.\eea The last equation provides inductively the
corestriction maps of a cohomomorphism ${\cal G}$. One can check
that ${\cal G}$ not only verifies ${\cal FG}={\cal I}$, but also
${\cal GF}={\cal I}$.\medskip

2. Take a $\op{Lod}_{\infty}$ structure $Q$ on $V$ and a
cohomomorphism ${\cal F}:(T(\da V),\zD)\to(T(\da V),\zD)$, such
that ${\cal F}_1=\op{id}$. Since $(f\otimes g)\circ (h\otimes
k)=(-1)^{\la g,h\ra}(f\circ h)\otimes (g\circ k)$, with
self-explaining notations, it is easily seen that ${\cal F}Q{\cal
F}^{-1}$, where ${\cal F}^{-1}$ is the inverse cohomomorphism
${\cal G}$ given by Item 1, is a weight $e_1$ codifferential of
$T(\da V)$, i.e. a $\op{Lod}_{\infty}$ structure on $V$.
Eventually, ${\cal F}$ is obviously a $\op{Lod}_{\infty}$
morphism, and, in view of the assumption ${\cal F}_1=\op{id}$,
even a $\op{Lod}_{\infty}$ isomorphism.\medskip

3. Consider two $\op{Lod}_{\infty}$ algebras $(V,Q)$, $(V',Q')$
and a $\op{Lod}_{\infty}$ isomorphism ${\cal F}$, i.e. a
cohomomorphism, such that ${\cal F}_1$ is bijective and $Q'{\cal
F}={\cal F}Q$ ($\star$). It then follows from Item 1 that there is
an inverse cohomomorphism ${\cal F}^{-1}$, such that $({\cal
F}^{-1})_1=({\cal F}_1)^{-1}$, and from Equation ($\star$) that
$Q{\cal F}^{-1}={\cal F}^{-1}Q'$.\end{proof}

The following key-theorem generalizes the last item of Proposition
\ref{InverMorph}.

\begin{theo}\label{Theo:QuaiInverse} If $f:(V,\pi)\raa (V',\pi')$
is a $\op{Lod}_{\infty}$ quasi-isomorphism, it admits a {\em
quasi-inverse}, i.e. there exists a $\op{Lod}_{\infty}$
quasi-isomorphism $g:(V',\pi')\raa(V,\pi)$, which induces the
inverse isomorphism in cohomology, i.e.
$g_{1\,\sharp}=(f_{1\,\sharp})^{-1}$.\end{theo}

We prove this theorem, which does not hold true in the category of
DGLodAs, in the next section.

\section{Minimal model theorem for Loday infinity algebras}

\begin{defi}
A $\op{Lod}_{\infty}$ algebra $(V,\pi)$ is {\em minimal}, if
$\pi_1=0$. It is {\em contractible}, if $\pi_p=0$, for $p\geq2,$
and if in addition $H(V,\pi_1)=0$.
\end{defi}

\begin{theo}\label{Theo:DecompOfLodAlgebra}
Each $\op{Lod}_{\infty}$ algebra is $\op{Lod}_{\infty}$ isomorphic
to the direct sum of a minimal $\op{Lod}_{\infty}$ algebra and a
contractible $\op{Lod}_{\infty}$ algebra.
\end{theo}

Theorem \ref{Theo:DecompOfLodAlgebra} was proved for
$\op{A}_{\infty}$ (resp. $\op{L}_{\infty}$) algebras in \cite{Ka}
(resp. in \cite{K} and e.g. \cite{AMM}). In the sequel, we provide
a proof in the $\op{Lod}_{\infty}$ case.

\begin{proof} Let $(V,\pi)$ be a $\op{Lod}_{\infty}$ algebra. For any
$\za\in\Z^n$, denote by $Z^{\za}$ and $B^{\za}$ the trace on
$V^{\za}$ of the kernel and the image of $\pi_1$. Consider a
supplementary vector subspace $V_m^{\za}$ of $B^{\za}$ in
$Z^{\za}$ and a supplementary subspace $W^{\za}$ of $Z^{\za}$ in
$V^{\za}$. Let $Z,B,V_m,$ and $W$ be the corresponding graded
spaces. Then, the complex $(V,\pi_1)$ decomposes into the direct
sum of the complex $(V_m,0)$, with vanishing differential, and the
complex $(V_c:=B\oplus W,\pi_1)$, with trivial cohomology. It
follows that the sequence $f^{(1)}:=(\op{id},0,\ldots)$ is a
$\op{Lod}_{\infty}$ isomorphism from $(V,\pi)$ to the
$\op{Lod}_{\infty}$ algebra $L_1:=(V_m\oplus
V_c,0\oplus\pi_1,\pi_2,\pi_3,\ldots).$ We will transform
inductively the maps $\pi_p$, $p\geq2$, via $\op{Lod}_{\infty}$
isomorphisms, into mappings of the form $\pi^m_p\oplus 0$, such
that $\pi^m:=(0,\pi_2^m,\zp_3^m,\ldots)$ be a minimal
$\op{Lod}_{\infty}$ structure on $V_m$.

\begin{lem}\label{Lemma:HomotopyOperator} Consider
the operator $\delta:V\raa V$ that is defined, for any $v\in
V_m\oplus W$, by $\delta(v)=0$, and, for any $v\in B$, by
$\delta(v)=w$, where $w$ is the unique element $w\in W$, such that
$\zp_1(w)=v$. Let $P$ be the projection onto $V_m$ with respect to
the decomposition $V=V_m\oplus V_c$. Then, $\delta$ is a homotopy
operator between the complex endomorphisms $P$ and $\op{id}$ of
$(V,\zp_1)$, i.e. $\pi_1\delta+\delta\pi_1=\op{id}-P.$
\end{lem}

\begin{proof} Obvious.\end{proof}

Let us construct $\pi_2^m$. Consider a sequence
$f^{(2)}:=(\op{id},f_2,0,0,\ldots)$, where $f_2$ is a weight
$-e_1$ bilinear map on $V$. According to Item 2 of Proposition
\ref{InveLodInftIsom}, $f^{(2)}$ defines a $\op{Lod}_{\infty}$
isomorphism
$$L_1\raa (V_m\oplus
V_c,\pi_1^{(2)},\pi_2^{(2)},\pi_3^{(2)},\ldots)=:(V_m\oplus
V_c,\pi^{(2)}),$$ and $\zp^{(2)}$ is a $\op{Lod}_{\infty}$
structure on $V_m\oplus V_c$, if and only if $\pi^{(2)}\circ
f^{(2)}=f^{(2)}\circ\pi.$ In view of Equation (\ref{75}), this
condition implies that (take $p=1$) $
\zp_1^{(2)}=\zp_1=0\oplus\zp_1,$ that (take $p=2$), for $v_1\in
V^{v_1},v_2\in V$,
\be\label{101}\pi^{(2)}_2(v_1,v_2)=-\pi_1f_2(v_1,v_2)+\pi_2(v_1,v_2)-f_2(\pi_1v_1,v_2)-(-1)^{\la
e_1,v_1\ra}f_2(v_1,\pi_1v_2),\ee and (take $p\ge 3$) it provides
the $\zp_p^{(2)}$, $p\ge 3,$ in terms of $f_2$.\medskip

It suffices to find a weight $-e_1$ bilinear map $f_2$, such that
the resulting $\zp_2^{(2)}$ maps $V_m\times V_m$ to $V_m$ and
vanishes elsewhere. Indeed, the restriction $\zp_2^m$ of
$\zp_2^{(2)}$ to $V_m\times V_m$ is then a weight $0$ bilinear map
on $V_m$ and $\zp_2^{(2)}=\zp_2^m\oplus 0$. If we choose the
$\zp_p^{(2)}$, $p\ge 3$, given by $f_2$, intertwining condition
$\pi^{(2)}\circ f^{(2)}=f^{(2)}\circ\pi$ is satisfied and
$f^{(2)}$ is a $\op{Lod}_{\infty}$ isomorphism between $L_1$ and
$(V_m\oplus V_c,0\oplus\pi_1,\pi_2^m\oplus
0,\pi_3^{(2)},\zp_4^{(2)},\ldots).$ When continuing step by step,
we finally get a $\op{Lod}_{\infty}$ isomorphism $\ldots
f^{(3)}\circ f^{(2)}\circ f^{(1)}$ between $(V,\zp)$ and
$(V_m\oplus V_c,0\oplus\pi_1,\pi_2^m\oplus 0,\pi_3^m\oplus
0,\ldots).$ It eventually follows from the explicit form of the
stem bracket, see Example \ref{Ex:DirectSumLodayInfty} and
subsequent explanation, that, since $(0\oplus\pi_1,\pi_2^m\oplus
0,\pi_3^m\oplus 0,\ldots)$ verifies the $\op{Lod}_{\infty}$
structure condition on $V_m\oplus V_c$, the two terms of this
direct sum verify the same condition on $V_m$ and $V_c$
respectively.\medskip

Let us define $f_2$ as follows: \be\label{Defif_2}
f_2(v_1,v_2)=\left\{
\begin{array}{lll} \delta\pi_2(v_1,v_2)+P\pi_2(w,v_2),\; \text{if}\;
(v_1,v_2)\in B^{\za}\times Z^{\zb},\;
v_1=\pi_1(w),\\\delta\pi_2(v_1,v_2)+\frac{1}{2}P\pi_2(w,v_2), \;
\text{if}\; (v_1,v_2)\in B^{\za}\times W^{\zb}, \;
v_1=\pi_1(w),\\\delta\pi_2(v_1,v_2)+(-1)^{\la
e_1,\za\ra}P\pi_2(v_1,w'),\; \text{if}\; (v_1,v_2)\in
Z^{\za}\times
B^{\zb},\\\quad\quad\quad\quad\quad\quad\quad\quad\quad\quad\quad\quad\quad\quad\quad\quad\quad\quad\quad\quad\quad
v_2=\pi_1(w'),\\\delta\pi_2(v_1,v_2)+(-1)^{\la
e_1,\za\ra}\frac{1}{2}P\pi_2(v_1,w'),\; \text{if}\; (v_1,v_2)\in
W^{\za}\times B^{\zb},\\\quad\quad\quad\quad\quad\quad\quad\quad\quad\quad\quad\quad\quad\quad\quad\quad\quad\quad\quad\quad\quad v_2=\pi_1(w'),\\
\delta\pi_2(v_1,v_2),\;\text{otherwise}.\end{array}\right.\ee Map
$f_2$ is well-defined, i.e. the two definitions on $B^{\za}\times
B^{\zb}\subset (B^{\za}\times Z^{\zb})\cap(Z^{\za}\times B^{\zb})$
coincide, as the $\op{Lod}_{\infty}$ structure condition
(\ref{LodaInftCond}) implies
$$0=P\zp_1\zp_2(v_1,v_2)=P\zp_2(\zp_1v_1,v_2)+(-1)^{\la
e_1,\za\ra}P\zp_2(v_1,\zp_1v_2),$$ for any $v_1\in V^{\za},v_2\in
V$. Indeed, when writing this upshot for $w\in W^{\za-e_1}$ and
$w'$, see Equation (\ref{Defif_2}), we get the announced
result.\medskip

It remains to show that $\zp_2^{(2)}$ sends $V_m\times V_m$ to
$V_m$ and vanishes elsewhere.\medskip

If $(v_1,v_2)\in Z\times Z$, Equation (\ref{101}) and Lemma
\ref{Lemma:HomotopyOperator} yield
$\pi^{(2)}_2(v_1,v_2)=\delta\pi_1\pi_2(v_1,v_2)+P\pi_2(v_1,v_2).$
But, since the $\op{Lod}_{\infty}$ condition entails that
$\pi_1\pi_2(v_1,v_2)=0,$ we get
$\pi^{(2)}_2(v_1,v_2)=P\pi_2(v_1,v_2).$ Furthermore, for any
$(w,v_2)\in W\times Z$, Condition (\ref{LodaInftCond}) implies
that $P\pi_2(\pi_1w,v_2)=0$, whereas for any $(v_1,w')\in Z\times
W$, we obtain $P\pi_2(v_1,\pi_1w')=0.$ Therefore,
\bea\pi_2^{(2)}(v_1,v_2)=\left\{ \begin{array}{lll}
P\pi_2(v_1,v_2)=:\pi^m_2(v_1,v_2)\in V_m, \; \text{if}
(v_1,v_2)\in V_m\times V_m, \\ 0, \; \text{if}\; (v_1,v_2)\in
B\times V_m\; \text{or}\; (v_1,v_2)\in V_m\times
B\;\text{or}\;(v_1,v_2)\in B\times B.
\end{array}\right.\eea

If $(v_1,v_2)\in (W\times Z)\cup (Z\times W)\cup (W\times W)$,
Equation (\ref{101}), Lemma \ref{Lemma:HomotopyOperator}, and
Condition (\ref{LodaInftCond}) allow checking that
$\pi^{(2)}_2(v_1,v_2)=0$.\medskip

Hence,
$$\pi_2^{(2)}(v_1,v_2)=\left\{ \begin{array}{lll}
\pi^m_2(v_1,v_2)=P\pi_2(v_1,v_2)\in V_m, \; \text{if} (v_1,v_2)\in
V_m\times V_m, \\ 0, \; \text{otherwise},\end{array}\right.$$ so
that it suffices to continue by induction.\end{proof}

We are now prepared to prove Theorem
\ref{Theo:QuaiInverse}.\medskip

\begin{proof} Let $f:(V,\pi)\raa(V',\pi')$ be a $\op{Lod}_{\infty}$
quasi-isomorphism. According to the minimal model theorem, there
is a $\op{Lod}_{\infty}$ isomorphism $h$ (resp. $h'$) that
identifies the $\op{Lod}_{\infty}$ algebra $(V,\pi)$ (resp.
$(V',\pi')$) to a direct sum $(V_m\oplus V_c,\pi^m\oplus\pi^c)$
(resp. $(V'_m\oplus V'_c,\pi'^m\oplus\pi'^c)$). Furthermore, since
the inclusion ${\frak i}:=(i,0,0,\ldots):(V_m,\pi^m)\raa
(V_m\oplus V_c,\pi^m\oplus\pi^c)$ (resp. the projection ${\frak
p}:=(P',0,0,\ldots):(V'_m\oplus V'_c,\pi^{\prime
m}\oplus\pi^{\prime c})\raa (V'_m,\pi'^{m}))$ is a
$\op{Lod}_{\infty}$ quasi-isomorphism, the sequence $h^{{\frak
i}}:=h^{-1}\circ{\frak i}$ (resp. $h'^{\frak p}:={\frak p}\circ
h'$) is a $\op{Lod}_{\infty}$ quasi-isomorphism from $(V_m,\pi^m)$
(resp. $(V',\pi')$) to $(V,\pi)$ (resp. $(V'_m,\pi'^m)$).
Therefore, the map $f^m:= h'^{\frak p}\circ f\circ h^{\frak i}$ is
a $\op{Lod}_{\infty}$ quasi-isomorphism between $(V_m,\pi^m)$ and
$(V_m',\pi'^m)$. But, as $H(V_m,\pi^m_1)=V_m$ and
$H(V_m',\pi_1'^m)=V'_m$, the map $(f_1^m)_{\sharp}=f_1^m:V_m\raa
V'_m$ is an isomorphism, and so $f^m$ has a $\op{Lod}_{\infty}$
isomorphism inverse $(f^m)^{-1}$, with
$(f^m)^{-1}_1=(f_1^m)^{-1}$, see Proposition
\ref{InveLodInftIsom}. Consequently, the sequence $g:=h^{{\frak
i}}\circ (f^m)^{-1}\circ h'^{\frak p}$ is a $\op{Lod}_{\infty}$
quasi-isomorphism from $(V',\pi')$ to $(V,\pi)$. Moreover,
$g_{1\sharp}=(f_{1\sharp})^{-1}$. Indeed, observe first that
$g_1=i\circ(f_1^m)^{-1}\circ P'$ and $f_1^m=P'\circ f_1\circ i$.
For any $[v']\in H(V',\zp'_1)$, we thus get
$g_{1\sharp}[v']=[(f_1^m)^{-1}P'v']\in H(V,\zp_1)$. On the other
hand, $(f_{1\sharp})^{-1}[v']=:[v]\in H(V,\zp_1)$, hence
$f_1v=v'+\zp'_1{\frak v}'$, ${\frak v}'\in V'.$ It now suffices to
show that there is ${\frak v}\in V$, such that
$(f_1^m)^{-1}P'v'=v+\zp_1{\frak v}$, i.e.
$$P'v'=f_1^m(v+\zp_1{\frak v})=P'f_1v+P'f_1\zp_1{\frak
v}=P'(v'+\zp'_1{\frak v}')+P'\zp'_1f_1{\frak v}.$$ This condition
is satisfied for any ${\frak v}\in V$, since $P'\zp'_1=0$.
\end{proof}

In view of the preceding proof, we have the following

\begin{cor}
Each $\op{Lod}_{\infty}$ algebra is $\op{Lod}_{\infty}$
quasi-isomorphic to a minimal one.
\end{cor}

\section{Graded and strongly homotopy algebra cohomologies}

\subsection{Graded Loday and Chevalley-Eilenberg
cohomologies}\label{SubSect:GradLodCoho}

Let $\zp\in\op{Lod}(V)$ be a $\Z^n$-graded Loday structure on $V$.
As $\zp$ is canonical for the $\Z^{n+1}$--GLA $(M_r(V),
[-,-]^{\otimes})$, see Theorem \ref{lodaycodifferential}, it is
clear that the cohomology of the induced DGLA, with differential
$\p_{\zp}=[\zp,-]^{\otimes}$, should roughly be the cohomology of
the considered Loday algebra.

\begin{prop}\label{GradLodaCoho}
The graded Loday cohomology operator $\p_{\zp}$ of a Loday
structure $\zp=\{-,-\}$ on a vector space $V$, reads, for any
$B\in M^{(B,b)}(V)$, $b\ge -1$, and any homogeneous $v_1,\ldots,
v_{b+2}\in V$,
\bean\label{lodaycohmologycoboundary}\nonumber&&(\p_{\pi}B)(v_1,\ldots,v_{b+2})=(-1)^{b+1}\{B(v_1,\ldots,v_{b+1}),v_{b+2}\}-
\\&&\sum_{i=1}^{b+1}(-1)^{i-1}(-1)^{\la B+v_1+\ldots+v_{i-1},v_{i}\ra}\nonumber
\{v_i,B(v_{1},\ldots,\widehat{v_i},\ldots,v_{b+1},v_{b+2})\}\\&& +
\sum_{i=1}^{b+1}\sum_{j=1}^i(-1)^{j+1}(-1)^{\la
v_j,v_{j+1}+\ldots+v_{i}\ra}\\&&
\nonumber\;\;\;\;\;\;\;\;\;\;\;\;\;\;\;\;\;\;\;\;B(v_{1},\ldots,\widehat{v}_j,\ldots,v_{i},
\{v_{j},v_{i+1}\},v_{i+2},\ldots,v_{b+2}).\eean
\end{prop}

\begin{proof} For $b\ge 0$, the explicit form of $\p_{\zp}$
is a consequence of Theorem \ref{ExplLodaBrac}. Equation
(\ref{lodaycohmologycoboundary}) suggests extending $\p_{\pi}$ to
$M^{-1}(V)=V$ by $(\p_{\pi}v)(w):=\pi(v,w)=\{v,w\},$ for any
$v,w\in V$. The extended operator $\p_{\zp}$ is a cohomology
operator on $M(V)$, since \bea
\p_{\pi}(\p_{\pi}v)(v_1,v_2)&=&-\pi(\zp(v,v_1),v_2)\\&&\quad-(-1)^{\la
v,v_1\ra}\pi(v_1,\zp(v,v_2))+\zp(v,\pi(v_1,v_2))=0,\eea  in view
of the Jacobi identity.\end{proof}

\begin{defi} The {\em graded Loday cohomology} of a $\Z^n$-graded Loday
algebra $(V,\zp)$ is the cohomology of the complex
$(M(V),\p_{\zp})$, where $\p_{\zp}$ is given by Equation
(\ref{lodaycohmologycoboundary}).\end{defi}

\begin{rem} In the non-graded case, Operator (\ref{lodaycohmologycoboundary})
coincides with the (non-graded) Loday coboundary operator, see
\cite{DT}, and in the antisymmetric situation, it is (the opposite
of) the graded Chevalley-Eilenberg differential, see \cite{LMS}.
\end{rem}

\subsection{Graded Poisson and Jacobi cohomologies}\label{sectionThe algebraic(Jacobi)Schoutenbracket}

In the following we prove that the stem bracket $[-,-]^{\otimes}$
not only restricts to the Nijenhuis-Richardson bracket, see Remark
\ref{Rem:StemToNR}, but also to the Grabowski-Marmo bracket, which
was defined in \cite{GM2}, and in particular to the
Schouten-Jacobi and Schouten brackets.

\subsubsection{Definition}

It is well-known that the triplet made up by the space of
multivector fields of a manifold, the wedge product, and the
Schouten-Nijenhuis bracket, is a graded Poisson algebra of weight
$\za=-1$. Let us recall that, more generally, a {\em graded
Poisson algebra} of weight $\za\in\Z^n$ is an associative
$\Z^n$-graded commutative algebra ${\cal A}$, endowed with a
bilinear bracket $\{-,-\}$ of weight $\za$, which is {\em
$\za$-graded antisymmetric}, i.e. verifies $\{u,v\}=-(-1)^{\la
u+\za, v+\za\ra}\{v,u\},$ and satisfies the graded Jacobi
identity, as well as the graded Leibniz rule
\be\{u,vw\}=\{u,v\}w+(-1)^{\la
u+\za,v\ra}v\{u,w\},\label{GradJacoIden}\ee for any homogeneous
$u,v,w\in\A$. The concept of {\em graded Jacobi algebra} of weight
$\za$ is defined similarly, except that the associative algebra
must have a unit $1$ and that Condition (\ref{GradJacoIden}) is
replaced by \be \{u,vw\}=\{u,v\}w+(-1)^{\la
u+\za,v\ra}v\{u,w\}-\{u,1\}vw.\label{GeneGradJacoIden}\ee Note
that we cannot confine our study to graded Jacobi (resp. Poisson)
structures of weight $\za=0,$ as weight $\za$ does not disappear
via an $\za$-shift in the grading of $V$.\medskip

Equation (\ref{GeneGradJacoIden}) (resp. Equation
(\ref{GradJacoIden})) means that a graded Jacobi (resp. Poisson)
structure of weight $\za$ is an $\za$-antisymmetric graded first
order bidifferential operator (resp. graded biderivation) of $\A$.
It follows that the appropriate graded Jacobi (resp. Poisson)
cochain space is the space $\op{Diff}_1({\cal A})$ (resp.
$\op{Der}({\cal A})$) of $\za$-antisymmetric graded first order
polydifferential operators (resp. graded polyderivations) of
$\A.$\medskip

In \cite{GM2}, the authors investigated these operators using a
variant of Krasil'shchik's calculus, see \cite{Kr2}, which is
based upon a particular bidegree of the spaces $M(\A)$, $A(\A)$,
or $\op{Diff}_1(\A)$ (resp. $\op{Der}(\A)$). More precisely, let
us fix $\za\in\Z^{n}$ (to avoid overcrowded notations, the
dependence on $\za$ of the below constructed objects will not be
specified). We then denote by ${\frak M}(\A)$ the usual vector
space $M(\A)$ with the $\Z^{n+1}$-gradation \be{\frak M}(\A)
=\bigoplus_{(A,a)\in\Z^n\times\Z}{\frak M}^{(A+\za
a,a)}(\A),\label{MultSpecBide}\ee where the space ${\frak
M}^{(A+\za a,a)}(\A)$ of bidegree $(A+\za a,a)$ vanishes for
$a<-1$, coincides with $\A^A$ for $a=-1$, and is, for $a\ge 0$,
the space of $(a+1)$-linear maps on ${\cal A}$ that have weight
$A$. The subspaces of $\za$-antisymmetric multilinear mappings
will be denoted by ${\frak A}^{(A+\za a,a)}(\A)$. Moreover, the
subspaces $\op{\frak Diff}_1^{(*,a)}(\A)$ (resp. $\op{\frak
Der}^{(*,a)}(\A)$) and the associative graded commutative
dot-product ``$\cdot$'' on these subspaces are described
inductively w.r.t. $a$, see \cite{GM2} and \cite{Kr2}. We then get
the following

\begin{prop} Let ${\cal A}$ be an associative
$\Z^n$-graded commutative unital algebra and let $\za\in\Z^n$.
With respect to the bidegree (\ref{MultSpecBide}), the pair
$(\op{\frak Diff}_1(\A),\cdot)$ is an associative
$\Z^{n+1}$-graded commutative unital algebra, $\cdot$ has weight
$(\za,1)$, and the pair $(\op{\frak Der}(\A),\cdot)$ is a
$\Z^{n+1}$-graded subalgebra.
\end{prop}

As the objective is to define the graded Jacobi (resp. Poisson)
cohomology, our task is to find on the cochain algebra $(\op{\frak
Diff}_1(\A),\cdot)$ (resp. $(\op{\frak Der}(\A),\cdot)$) a
$\Z^{n+1}$-graded Lie structure (of weight $(0,0)\in\Z^n\times\Z$)
or even a graded Jacobi (resp. Poisson) structure. If such a
bracket $[-,-]^{\op{GM}}$ exists and satisfies
\be[v,w]^{\op{GM}}=0,
[A,v]^{\op{GM}}=(-1)^{1-a}A(v),\label{DefGM0}\ee it necessarily
also verifies \be [v,A]^{\op{GM}}=-(-1)^{\la v-\za,A+\za
a\ra+a}[A,v]^{\op{GM}}\label{DefiGM1}\ee and
\be\label{DefiGM2}\begin{array}{lll}[A,B]^{\op{GM}}(v)&=&(-1)^{1-a-b}[[A,B]^{\op{GM}},v]^{\op{GM}}\\&=&(-1)^{1-a-b}
[A,[B,v]^{\op{GM}}]^{\op{GM}}+(-1)^{1+\la B+\za b,v-\za\ra
-a}[[A,v]^{\op{GM}},B]^{\op{GM}}\\&=&
(-1)^{a}[A,B(v)]^{\op{GM}}+(-1)^{\la B+\za
b,v-\za\ra}[A(v),B]^{\op{GM}},\end{array}\ee where $v,w$ (resp.
$A,B$) are homogeneous elements of ${\cal A}\subset\op{\frak
Diff}_1({\cal A})$ (resp. of $\op{\frak Diff}_1({\cal A})$). The
bracket that is defined inductively by Equations (\ref{DefGM0}),
(\ref{DefiGM1}), and (\ref{DefiGM2}) actually fits, see \cite{GM2}
and \cite{Kr2}:

\begin{theo} If ${\cal A}$ is an associative
$\Z^n$-graded commutative unital algebra and $\za\in\Z^n$, there
is a unique $\Z^{n+1}$-graded Jacobi bracket $[-,-]^{\op{GM}}$ of
degree $(0,0)\in\Z^n\times\Z$ on the associative $\Z^{n+1}$-graded
commutative algebra $(\op{\frak Diff}_1(\A),\cdot)$ that verifies
$[A,v]^{\op{GM}}=(-1)^{1-a}A(v)$, $A\in\op{\frak
Diff}_1^{(*,a)}({\cal A})$, $v\in{\cal A}$. Moreover, $(\op{\frak
Der}(\A),[-,-]^{\op{GM}},\cdot)$ is a $\Z^{n+1}$-graded Poisson
algebra. Furthermore, the graded Jacobi $($resp. Poisson$)$
structures of weight $\za$ on $\A$, are exactly the canonical
elements of degree $(2\za,1)$ of this graded Jacobi $($resp.
Poisson$)$ algebra.\label{CanoGradJaco}\end{theo}

This result immediately leads to the

\begin{defi} If $\zp$ denotes a graded Jacobi $($resp. Poisson$)$
structure of weight $\za\in\Z^n$ on the (usual) algebra $\A$, we
refer to the cohomology of the DGLA $(\op{\frak
Diff}_1(\A),[-,-]^{\op{GM}},[\zp,-]^{\op{GM}})$ (resp. $(\op{\frak
Der}(\A),[-,-]^{\op{GM}},[\zp,-]^{\op{GM}})$) as the {\em graded
Jacobi $($resp. Poisson$)$ cohomology} of $(\A,\zp)$.
\end{defi}

Of course, in the (ungraded) geometric case $\A=\Ci(M)$ (usual
notations), the graded Jacobi (resp. Poisson) cohomology coincides
with the standard Lichnerowicz-Jacobi (resp. Lichnerowicz-Poisson)
cohomology, see \cite{Lic} (resp. see \cite{AL}).

\subsubsection{Link with the stem bracket}

\begin{prop}\label{GMbrachetProp1} For any fixed $\za\in\Z^n$, any $A\in\op{\frak Diff}_1^{(A+\za a,a)}(\A)$
and any $B\in\op{\frak Diff}_1^{(B+\za b,b)}(\A)$, the
Grabowski-Marmo bracket is given by \be [A,B]^{\op{GM}}=A\square
B-(-1)^{\la A+\za a,B+\za b\ra+ab} B\square
A,\label{GMbracket1}\ee where $A\square B$ is defined inductively
by $$v\square w:=0,\quad A\square v:=(-1)^{1-a}A(v),\; a\ge
0,\quad v\square A:=0,\; a\ge 0,$$ \be (A\square
B)(v):=(-1)^{a}A\square B(v)+(-1)^{\la B+\za
b,v-\za\ra}A(v)\square B,\; a,b\ge 0,\label{DefiSqua}\ee for any
$v\in\A^v$, $w\in\A$. Further, the explicit form of the
square-product is \be\begin{array}{ll} (A\square
B)\!\!\!\!&(v_1,\ldots,v_{a+b+1})
\\&=(-1)^{1+ab}\sum_{\begin{array}{cc}I\cup J=N^{(a+b+1)}\\\vert I\vert=b+1,\vert
J\vert=a\end{array}}(-1)^{(I;J)}\ze_{\da
V}(I;J)A(B(V_I),V_J),\end{array}\label{SquaExplForm}\ee where
$v_i\in\A^{v_i}$, $\da V=\da v_1\otimes\ldots\otimes\da
v_{a+b+1}$, and $\da v_i\in \A^{v_i-\za}$.
\end{prop}

\begin{proof} The proof of Equation (\ref{GMbracket1}) is by induction on
$a+b$. If $a+b\le -1$, the claim is obviously true, see Equations
(\ref{DefGM0}) and (\ref{DefiGM1}). Equations (\ref{DefiGM2}) and
(\ref{DefiSqua}) allow checking the rest. In order to determine
the explicit form (\ref{SquaExplForm}), we proceed similarly. The
announced upshot is clear for $a=-1$ or $b=-1$, hence, in
particular for $a+b\le -1$. Let us now assume that it is valid for
$a+b\le k-1$, $k\ge 0,$ and examine the case $a+b=k.$ As already
pointed out, if $a=-1$ or $b=-1$, our conjecture is verified.
Otherwise, $a,b\ge 0$, and, if we set $V''=v_2\otimes\ldots\otimes
v_{a+b+1}$, if we use Equation (\ref{DefiSqua}) and the induction
assumption, and if we exchange $v_1$ and $B(V''_I)$ in the second
term of the result produced in this way, we obtain
$$\begin{array}{lll}&(A\square B)(v_1,v_2,\ldots,v_{a+b+1})&\\
=&(-1)^{1+ab}[\,\sum_{\begin{array}{l}I\cup J=\{2,\ldots,a+b+1\}\\
\vert I\vert=b,\vert J\vert=a\end{array}}(-1)^{(I;J)}\ze_{\da
V''}(I;J)A(B(v_1,V''_I),V''_J)&\\&+(-1)^{1+\la V''_I-\za
(b+1),v_1-\za\ra+b}&\\&\sum_{\begin{array}{l}I\cup
J=\{2,\ldots,a+b+1\}\\\vert I\vert=b+1,\vert
J\vert=a-1\end{array}}(-1)^{(I;J)}\ze_{\da
V''}(I;J)A(B(V''_I),v_1,V''_J)].\end{array}$$ In the final upshot
for which we look, ${\cal I}$ and ${\cal J}$ are unshuffles that
form a partition of $N^{(a+b+1)}$, so that $1$ is either the first
element of ${\cal I}$ or of ${\cal J}$. The first term
$\sum\ldots$ (resp. second term $(-1)^{\ldots}\sum\ldots$) of the
RHS of the last equation, corresponds exactly to the first (resp.
second) possibility. Indeed,
$$(-1)^{(I;J)}\ze_{\da V''}(I;J)=(-1)^{(\cal I;J)}\ze_{\da
V}({\cal I;J})$$ (resp. $$(-1)^{1+\la V''_I-\za
(b+1),v_1-\za\ra+b} (-1)^{(I;J)}\ze_{\da V''}(I;J)=(-1)^{(\cal
I;J)}\ze_{\da V}({\cal I;J})),$$ where $V=v_1\otimes\ldots\otimes
v_{a+b+1}$.\end{proof}

In the following, we denote by $\da\A$ the vector space $\A$
endowed with the gradation $(\da\A)^{\zg}=\A^{\zg+\za}$. Set now
\be{\frak A}(\A)=\bigoplus_{(A,a)\in\Z^n\times\Z}{\frak A}^{(A+\za
a,a)}(\A) \;\left(\mbox{resp. }\;
A(\da\A)=\bigoplus_{(A,a)\in\Z^n\times\Z}A^{(A+\za
a,a)}(\da\A)\right),\label{TwoGardDefi}\ee

\vspace{2mm}\noindent where ${\frak A}^{(A+\za a,a)}(\A)$ (resp.
$A^{(A+\za a,a)}(\da\A)$) is the space of multilinear maps
$A:\A^{\times(a+1)}\to\A$ (resp. $\tilde
A:(\da\A)^{\times(a+1)}\to\da\A$) that have weight $A$ (resp.
$A+\za a$), and that are further $\za$-graded antisymmetric (resp.
graded antisymmetric). Obviously, the map $\sim:A\to \tilde A$,
with $\tilde A$ defined by $\tilde A(\tilde v_1,\ldots,\tilde
v_{a+1})=\da A(\ua\tilde v_1,\ldots,\ua\tilde v_{a+1}),$ is a
$\Z^{n+1}$-graded vector space isomorphism, the inverse of which
is $\backsim:\tilde A\to A$, where $A(v_1,\ldots, v_{a+1})=\ua
\tilde A(\da v_1,\ldots, \da v_{a+1}).$ The isomorphism $\sim$
pulls of course the Nijenhuis-Richardson graded Lie bracket
$[-,-]^{\op{NR}}$ on the usual space $A(\da\A)$, associated with
the $\Z^n$-graded vector space $\da\A$, back to a graded Lie
bracket
$$\backsim\!\![\sim-,\sim-]^{\op{NR}}=:\,\sim^*\!\![-,-]^{\op{NR}}
=\,\sim^*\!\![-,-]^{\otimes}\vert_{A(\da\A)}$$ on $\op{\frak
A}(\A)$.

\begin{prop} The Grabowski-Marmo bracket $[-,-]^{\op{GM}}$ is the restriction of
$\sim^*\!\![-,-]^{\otimes}\vert_{A(\da\A)}$ to $\op{\frak
Diff}_1(\A)$.
\end{prop}

\begin{proof} The explicit form of the restriction of the stem
bracket, see Theorem \ref{ExplLodaBrac}, to skew-symmetric
mappings, i.e. the explicit expression of the Nijenhuis-Richardson
bracket, as well as the forms (\ref{SquaExplForm}) and
(\ref{GMbracket1}) of the Grabowski-Marmo bracket imply that, for
$A\in \op{\frak Diff}_1^{(A+\za a,a)}(\A)$ and $B\in \op{\frak
Diff}_1^{(B+\za b,b)}(\A)$, the preceding pullback reads
$$\lp\sim^*\!\![-,-]^{\otimes}\vert_{A(\da\A)}\rp(A,B)=\backsim i_{\tilde A}\tilde B-(-1)^{\la A+\za a,B+\za
b\ra+ab}\backsim i_{\tilde B}\tilde A,$$ where $\backsim i_{\tilde
B}\tilde A=(-1)^{1+\la A+\za a,B+\za b\ra+ab}A\square B,$ so that
the result follows.\end{proof}

Finally:

\begin{cor} The stem bracket is (up to reading through a canonical isomorphism)
a graded Lie bracket on the spaces of graded Loday, graded Lie,
graded Poisson, and graded Jacobi cochains, for which the
corresponding algebraic structures are canonical elements, and
that thus encodes the graded cohomologies of all these structures.
\end{cor}

\subsection{Strongly homotopy and graded $\mathbf{p}$-ary Loday cohomologies}

Since $\op{Lod}^Q_{\infty}$ structures on a $\Z^n$-graded vector
space $V$, $Q\in\Z^n$, $\la Q,Q\ra$ odd, are the degree $Q$
canonical elements $\zp$ of the $\Z^n$-graded Lie algebra
$(C(V),[-,-]^{\bar{\otimes}})$, we have the natural

\begin{defi} The {\em cohomology of a Loday infinity algebra} $(V,\zp)$ is
the cohomology of the DGLA
$(C(V),[-,-]^{\bar{\otimes}},[\zp,-]^{\bar{\otimes}})$, where the
coboundary operator is, for any $\zr\in C^{\zr}(V),$ explicitly
given by
$$\begin{array}{ll}[\zp,\zr]^{\bar{\otimes}}=\sum_{q\ge 1}\sum_{s+t=q+1}(-1)^{1+(s-1)\la e_1,\zr
\ra}[\zp_s,\zr_t]^{\otimes}\end{array}$$ (and Equations
(\ref{bracketloday}) and (\ref{MLodaycomposition})).
\end{defi}

1. We first examine the case $\zp=\zp_p$, $p\in\{2,4,\ldots\}$. If
the odd degree $Q\in\Z^n$ of $\zp$ is chosen to be equal to
$(p-1)e_1$, we have \be\zp=\zp_p\in M^{(0,p-1)}(V),
[\zp_p,\zp_p]^{\otimes}=0, p\; \mbox{even}.\label{nAryLodaAlge}\ee

Essentially two $p$-ary extensions of the Jacobi identity were
investigated during the last decades. If $[-,-,\ldots,-]$ denotes
a $p$-linear bracket on $V$, the first is the generalization,
which requires that the adjoint action
$[v_1,v_2,\ldots,v_{p-1},-]$ be a derivation for the $p$-ary
bracket $[w_1,w_2,\ldots,w_p]$, see e.g. \cite{Fil}, and leads to
{\em Nambu-Lie} or, in the nonantisymmetric context, to {\em
Nambu-Loday} structures, see \cite{Nam}. The second was suggested
by P. Michor and A. Vinogradov, see \cite{MichVino}, and further
studied in \cite{VinoVino} and in \cite{VinoVino2}. We refer to
this last $p$-ary extension as {\em $p$-ary Lie} or {\em $p$-ary
Loday} structure. A {\em $p$-ary $\Z^n$-graded Lie structure} on a
$\Z^n$-graded vector space $V$ is a map $P_p\in A^{(0,p-1)}(V)$,
such that $i_{P_p}P_p=0$ (where $i$ is the interior product used
in the definition of the graded Nijenhuis-Richardson bracket) or,
if $p$ is even, equivalently, such that $[P_p,P_p]^{\op{NR}}=0.$
Moreover, the cohomology of a $p$-ary Lie algebra $(V,P_p)$ is the
cohomology of the DGLA $(A(V),[-,-]^{\op{NR}},[P_p,-]^{\op{NR}})$,
see \cite{MichVino}.\medskip

Analogously, a map $\zp_p\in M^{(0,p-1)}(V)$,
$p\in\{2,4,\ldots\}$, that verifies $[\zp_p,\zp_p]^{\otimes}=0,$
is a {\em $p$-ary $\Z^n$-graded Loday structure} on $V$. The
cohomology of such an algebra $(V,\zp_p)$ should be defined (and
was defined in the case $p=2$ roughly) as the cohomology of the
DGLA $(M_r(V),[-,-]^{\otimes},[\zp_p,-]^{\otimes})$.\medskip

In the following, we explain why

\begin{rem} The cohomology space of $\zp=\zp_p$, see Equation
(\ref{nAryLodaAlge}), viewed as degree $(p-1)e_1$ strongly
homotopy Loday structure, coincides with the preceding cohomology
space of $\zp=\zp_p$, viewed as $p$-ary $\Z^n$-graded Loday
structure.\end{rem}

Let us first mention that for noninfinity algebras, it is
conventional to substitute in cochain space $C(V)$ a direct sum
for the direct product, so that
$$C(V)=\bigoplus_{s\in\N^*}\bigoplus_{Q\in\Z^n}M^{(Q+(1-s)e_1,s-1)}(V)=\bigoplus_{s\in\N^*}M^{s-1}(V)=M_r(V).$$
Hence, $C(V)$ carries a bigrading that is shifted with respect to
the usual bigradation
$$M_r(V)=\bigoplus_{s\in\N^*}\bigoplus_{Q\in\Z^n}M^{(Q,s-1)}(V)$$ of
$M_r(V).$ We emphasize that in the first bigrading the space of
bidegree $(Q,s)$ is, unlike in Section \ref{sectionThe
algebraic(Jacobi)Schoutenbracket}, but just as above that section,
the space $M^{(Q+(1-s)e_1,s-1)}(V)$ -- whereas no ambiguity is
possible for the second grading.

It is now easy to see that the cohomology spaces of
$(C(V),[-,-]^{\bar{\otimes}},[\zp_p,-]^{\bar{\otimes}})$ and
$(M_r(V),[-,-]^{\otimes},$ $[\zp_p,-]^{\otimes})$ coincide.
Indeed, if $\zp$ is a $\op{Lod}^Q_{\infty}$ structure on $V$ and
if $\zr_t\in M^{(\zr+(1-t)e_1,t-1)}(V)$, we have \be
[\zp,\zr_t]^{\bar{\otimes}}=\sum_{q\ge 1}(-1)^{1+(q-t)\la e_1,\zr
\ra}[\zp_{q-t+1},\zr_t]^{\otimes}\in\prod_{q\ge
1}M^{(Q+\zr+(1-q)e_1,q-1)}(V),\label{LodaInfiCoboActi}\ee so that,
in the case $\zp=\zp_p$, $Q=(p-1)e_1$, where necessarily
$q=t+p-1$, the weight of cohomology operator
$[\zp_p,-]^{\bar{\otimes}}$ with respect to the first mentioned
bigrading is $((p-1)e_1,p-1)$. On the other hand, since $\zp_p\in
M^{(0,p-1)}(V)$, it is clear that the weight of operator
$[\zp_p,-]^{\otimes}$ with respect to the second bigradation is
$(0,p-1)$. It follows that the cohomology spaces of
$(C(V),[\zp_p,-]^{\bar{\otimes}})$ and
$(M_r(V),[\zp_p,-]^{\otimes})$, say $\bar{H}$ and $H$, are both
$\Z^{n+1}$-graded. Space $\bar{H}^{(\zr+(t-1)e_1,t-1)}$ (resp.
$H^{(\zr,t-1)}$) is encoded in the cocycle equation
$[\zp_p,\zr_t]^{\bar{\otimes}}=0$, i.e.
$[\zp_p,\zr_t]^{\otimes}=0,$ where $\zr_t\in M^{(\zr,t-1)}(V)$
(resp. $[\zp_p,\zr_t]^{\otimes}=0,$ $\zr_t\in M^{(\zr,t-1)}(V)),$
and in the coboundary equation, which reads
$\zr_t=[\zp_p,\zt]^{\bar{\otimes}}=[\zp_p,\pm\zt]^{\otimes},$ with
$\zt\in M^{(\zr,t-p)}(V)$ (resp. $\zr_t=[\zp_p,\zt]^{\otimes},$
$\zt\in M^{(\zr,t-p)}(V)),$ for any cocycle $\zr_t$. This
observation entails that
$$\bar{H}^{(\zr+(t-1)e_1,t-1)}=H^{(\zr,t-1)}.$$ Hence, the cohomology
spaces $\bar{H}$ and $H$ coincide, but their natural GLA
structures, which are induced by $[-,-]^{\bar{\otimes}}$ and
$[-,-]^{\otimes}$, are $\Z^n$- and $\Z^{n+1}$-graded
respectively.\medskip

For $p=2$, we of course recover the abovedefined cohomology space
of a graded Loday algebra.\medskip

2. If $\zp$ is not a sequence of all but one vanishing elements,
Equation (\ref{LodaInfiCoboActi}) implies that the Loday infinity
coboundary operator $[\zp,-]^{\bar{\otimes}}$ maps
$M^{(\zr+(1-t)e_1,t-1)}(V)$ into $\prod_{q\ge
1}M^{(Q+\zr+(1-q)e_1,q-1)}(V)$. Hence, the Loday infinity
cohomology is not $\Z$-graded. Nevertheless, if we consider the
(decreasing) filtration
$$C_k(V)=\bigoplus_{R\in\Z^n}\prod_{s\ge k}M^{(R+(1-s)e_1,s-1)}(V),\;k\ge 1,$$
and if $\zr=\sum_R\sum_{t\ge k}\zr_{R,t}\in C_k(V)$, we get
$$\begin{array}{ll}[\zp,\zr]^{\bar{\otimes}}=&\sum_R\sum_{q\ge 1}\sum_{s+t=q+1}(-1)^{1+(s-1)\la
e_1,R\ra}[\zp_{s},\zr_{R,t}]^{\otimes}\\\\&\quad\in\bigoplus_{R\in\Z^n}\prod_{q\ge
k}M^{(Q+R+(1-q)e_1,q-1)}(V).\end{array}$$ Indeed, as $k\le t$, we
have $q\le k-1\Rightarrow q\le t-1\Leftrightarrow s=q-t+1\le 0$,
so that the sum over $s,t$ vanishes for these $q$. The observation
yields $$[\zp,C_k(V)]^{\bar{\otimes}}\subset C_k(V).$$

Eventually,

\begin{rem} The pair $(C(V),[\zp,-]^{\bar{\otimes}})$ is a
differential filtered space and admits a spectral sequence.\\
\end{rem}

3. Let ${\cal C}(V)$ be the $\Z^{n}$-graded vector subspace
$${\cal C}(V)=\bigoplus_{Q\in\Z^n}{\cal
C}^Q(V)=\bigoplus_{Q\in\Z^n}\prod_{s\ge 1}
A^{(Q+(1-s)e_1,s-1)}(V)$$ of $C(V)$. Note that ${\cal C}(V)$ is
closed for the bracket $[-,-]^{\bar\ot}$ and that Lie infinity
structures are canonical elements of the graded Lie subalgebra
$({\cal C}(V),[-,-]^{\bar\ot})$. Hence, the {\em cohomology of a
Lie infinity algebra} $(V,\zp)$ is the cohomology of the DGLA
$({\cal C}(V),$ $[-,-]^{\bar{\otimes}},$
$[\zp,-]^{\bar{\otimes}})$. If the $\op{L}_{\infty}$ structure has
the degree $e_1$, the restriction of the coboundary operator
$[\zp,-]^{\bar\ot}$ to ${\cal C}(V)$ coincides with the
differential studied in \cite{Pen} and used in \cite{FP}.

\section{Acknowledgements} The authors thank S. Gutt for her interest in this work, fruitful discussions, and enriching
comments.

\end{document}